 \documentclass[preprint,12pt]{elsarticle}
\usepackage{amssymb}
\usepackage{graphicx}
\usepackage{subcaption}
\usepackage{amsmath}

\newcounter{theorem}
\newtheorem{theorem}{Theorem}
\newcounter{lemma}
\newtheorem{lemma}{Lemma}
\newcounter{corollary}
\newtheorem{corollary}{Corollary}
\newcounter{example}
\newtheorem{example}{Example}
\newcounter{proposition}
\newtheorem{proposition}{Proposition}

\begin{document}

\begin{frontmatter}

% Title, authors and addresses

% use the thanksref command within \title, \author or \address for footnotes;
% use the corauthref command within \author for corresponding author footnotes;
% use the ead command for the email address,
% and the form \ead[url] for the home page:
% \title{Title\thanksref{label1}}
% \thanks[label1]{}
% \author{Name\corauthref{cor1}\thanksref{label2}}
% \ead{email address}
% \ead[url]{home page}
% \thanks[label2]{}
% \corauth[cor1]{}
% \address{Address\thanksref{label3}}
% \thanks[label3]{}

\title{Curvature of planar aesthetic curves}

% use optional labels to link authors explicitly to addresses:
% \author[label1,label2]{}
% \address[label1]{}
% \address[label2]{}

\author{A. Cant\'on, L. Fern\'andez-Jambrina, M.J. V\'azquez-Gallo} \address{Departamento de Matem\'atica e
Inform\'atica Aplicadas a las Ingenier\'{\i}as Civil y Naval\\
Universidad Polit\'ecnica de Madrid\\
Avenida de la Memoria 4\\ E-28040-Madrid, Spain}

\begin{abstract}
% Text of abstract
In \cite{farin-a} Farin proposed a method for designing B\'ezier
curves with monotonic curvature and torsion.  Such curves are relevant
in design due to their aesthetic shape.  The method relies on applying
a matrix $M$ to the first edge of the control polygon of the curve in
order to obtain by iteration the remaining edges.  With this method,
sufficient conditions on the matrix $M$ are provided, which lead to
the definition of Class A curves, generalising a previous result by
Mineur et al.  \cite{mineur} for plane curves with $M$ being the
composition of a dilatation and a rotation.  However, Cao and Wang
\cite{cao} have shown counterexamples for such conditions.  In this
paper, we revisit Farin's idea of using the subdivision algorithm to
relate the curvature at every point of the curve to the curvature at
the initial point in order to produce a closed formula for the
curvature of planar curves in terms of the eigenvalues of the matrix
$M$ and the seed vector for the curve, the first edge of the control
polygon.  Moreover, we give new conditions in order to produce planar
curves with monotonic curvature.  The main difference is that we do
not require our conditions on the eigenvalues to be preserved under
subdivision of the curve.  This facilitates giving a unified derivation of the existing results and obtain more general results
in the planar case.
\end{abstract}

\end{frontmatter}

% main text
\section{Introduction}
\label{intro}

In many design applications, aesthetically pleasing curves and
surfaces are preferred.  The term ``class A surface" appeared when
Mercedes-Benz CAD/CAM systems' developers required a high quality
shape from the ``outside surface" parts, for which ``Au\ss enhaut" is
the German word, but they have been used in other industries, such as
naval architecture.

For a curve to be considered pleasing, its curvature and torsion
should vary monoto\-nical\-ly with respect to the parameter of the
curve.  According to \cite{farin}, a curve is \emph{fair} if its
curvature plot is continuous and consists of only a few monotone
pieces. There are several papers which deal with this issue for
B\'ezier and B-spline curves.  See \cite{sapidis, fairshape,
miura-review} for nice compilations on this issue.  In \cite{levien}
the issue of the \emph{best} spline curve is discussed.  One of the
possibilities are log-aesthetic curves and their generalisations
\cite{miura-review}.  Another approach deals with curves and surfaces
minimising functionals related to elasticity theory \cite{moreton,
moreton1}.  And another choice is restriction to a known family of
curves such as clothoids, as in civil engineering,
\cite{meek-clothoid, meek}, or logarythmic spirals.  In this paper we
focus on the class of curves known as ``Class A curves"
\cite{farin-a}.  In the context of regression functions, the
requirement of monotonicity and convexity has lead to the concept of
C-splines \cite{meyer}.

%
% \cite{andersson, beeker, farin-benz, hagen,
% hoschek-a, meier, miura, reyes-a, roulier}.  In
% \cite{miura-review} there is another review with special emphasis on
% log-aesthetic curves and their generalisation.

In \cite{mineur} it is defined a {\it typical curve}, a 2D B\'ezier
curve for which each edge of its control polygon is obtained after
rotating a given angle and changing the length of the previous edge by
a given factor.  The curve has monotonic curvature when a simple
condition involving the angle rotated and the scale factor is
satisfied.  Previous results on this issue may be found in
\cite{higashi}.

In \cite{farin-a} Farin discusses a 3D generalization of this concept
of typical curves and defines a special B\'ezier curve
by its degree, a vector $\mathbf{v}$ describing its first edge and a
square matrix $M$, in such a way that each edge is obtained from the
previous one by multiplying it by $M$.  Two conditions for the matrix
$M$ are given that are claimed to be sufficient to induce a Class A
B\'ezier curve with monotonic curvature and torsion.  These Class A
matrices are viewed as ``expansion" matrices that ``do not distort"
lengths ``too much".

These two conditions on a square matrix are supposed to imply two
pro\-per\-ties of the induced Class A curve which guarantee
monotonicity of curvature and torsion.  The first property is that
subdivision of a Class A B\'ezier curve provides a Class A B\'ezier
curve.  The second property requires that for a Class A B\'ezier curve the
curvature and torsion at the initial point are not smaller (greater) than
the ones at the endpoint.  In other words, the curvature and torsion
of a Class A curve are decreasing (increasing) functions of the
parameter of the curve.

However, in \cite{cao} it is claimed that Class A B\'ezier curves defined in
\cite{farin-a} are not invariant under subdivision based on a
counterexample in which the second condition of a Class A matrix (that of
 not distorting lengths too much) is supposed not to be invariant
under subdivision.  They also point out that the proof of this
fact in \cite{farin-a} is incomplete.  They propose new conditions
for a matrix to be called a Class A matrix in the case of symmetric
matrices and verify that these are sufficient conditions for
the induced Class A B\'ezier curve to have the two properties that
guarantee monotonicity of its curvature and torsion.  In the case of a
non-symmetric matrix, the authors state that their new conditions are
not sufficient to guarantee monotonic curvature and torsion.

Furthermore, in \cite{zhao} it is stated that the conditions defining Class A
matrices in \cite{farin-a} are incorrect using two counterexamples in
which Farin's conditions are supposed to hold for a curve, but its
curvature is not monotonic.  The issue of finding out necessary and sufficient
conditions to guarantee monotonicity of curvature is seen as future
work.

Finally, in \cite{hou} geometric sufficient conditions are provided for
monotonic curvature, but with no reference to Class A matrices.

In this paper we recover Farin's idea of using the subdivision
algorithm for relating the curvature at a point on the curve to the
curvature at the initial point, producing a new closed formula for the
curvature of planar curves in terms of the eigenvalues of the matrix
$M$ and the seed vector for the curve.  We also obtain new and simple
conditions on the eigenvalues of the matrix $M$, instead of its
singular values, for designing planar aesthetic curves.  In this we
generalise and depart from Farin's framework, since our conditions are
not preserved under subdivision and involve the seed vector, while
comprising previous results by Mineur et al (see \cite{mineur}) and Cao and Wang (see \cite{cao}).

Moreover, we discuss several claims present in \cite{farin-a},
\cite{cao} and \cite{zhao}.  We explicitly show that in \cite{farin-a}
neither the condition on the singular values of $M$ is the right one
for achieving a monotonic curvature nor it is preserved under
subdivision of the curve.  We note that the example in \cite{cao} for
a Class A matrix in which the condition involving singular values of
the matrix is not supposed to be invariant under subdivision is not a
strict counterexample for Farin's method.  Finally, we point out that
the first example of \cite{zhao} is not a counterexample for Farin's
method.

This paper is organised in the following way.  In
Section~\ref{baseautovectores} we disclose a new general formula for
the curvature of a planar curve in terms of the eigenvalues of the
matrix $M$ and the seed vector $\mathbf{w}$. This enables us to
provide new conditions for designing planar aesthetic curves in
Section~\ref{newconditions}.  Several examples are given.  Finally, in
Section~\ref{comparison} we discuss previous issues on Class A
curves in the literature.

% review some definitions
% and proofs presented in \cite{farin-a} and discuss their validity and
% relate them to the counterexamples in \cite{cao} and \cite{zhao}.  In
% Section~\ref{baseautovectores} we disclose a new formula for the
% curvature of Class A planar curves in terms of the eigenvalues of the
% matrix $M$ and the seed vector.
% is devoted to the
% validity of formula (10) in \cite{farin-a} and to the singular values
% of a matrix satisfying first condition of \cite{farin-a} and to the
% examples given in \cite{cao, zhao}.  We consider in
% Section~\ref{condition} the particular but important case of establishing conditions
% for a vector and a square matrix to induce a plane cubic B\'ezier
% curve with monotonic curvature.  A final section of conclusions is
% included.

\section{Curvature of planar B\'ezier curves generated by a matrix
\label{baseautovectores}}

In \cite{mineur}, Mineur et al. consider planar B\'ezier
curves for which each edge of the control polygon is obtained from the
previous one by the action of a square matrix $M$ that consists on a rotation and scaling, what they call \emph{typical curves}. They give
a relation between the rotation angle and the scaling factor in order to obtain B\'ezier curves with monotonic curvature. Farin (see \cite{farin-a}) generalizes this approach
to space curves and more general matrices $M$ with the goal to construct curves with monotonic curvature and torsion.

Following \cite{mineur} and \cite{farin-a} we study planar B\'ezier curves for which the edges of its control polygon are generated by the action of a
general square matrix $M$ on a vector $\mathbf{w}\neq \mathbf{0}$, the first edge of the control polygon. Concretely, given a $2\times 2$ matrix $M$ and a vector
$\mathbf{w} \neq \mathbf{0}$ in $\mathbb{R}^2$, we consider B\'ezier
curves of degree $n$ parameterised as
\[
\label{beziercurve} c(t)=\sum_{j=0}^n \, b_jB_j^n(t)	\qquad t\in[0,1],
\]
where $b_j$ are the control points of the curve, $B_j^n$ are the
Bernstein polynomials of degree $n$ and the edges of the control polygon
are given by
\begin{equation} b_{j+1}-b_j=M^j  \mathbf{w},\qquad j=0, \dots, n-1.
\end{equation}

Recall that the curvature of a parameterised planar curve is
\[
\kappa(t)=\dfrac{\det\bigl(c'(t),c''(t)\bigr)}{\Vert c'(t)\Vert^3}.
\]
Following \cite{farin}, we may write the expressions for the curvature at the endpoints of the
curve in terms of the edges of the control polygon,
%\begin{equation}\label{curvature0}\kappa(0)=2 \frac{n-1}{n}\frac{| D|}{\|\mathbf{w} \|^3}
%=\frac{n-1}{n}\frac{|\det(M\mathbf{w},\mathbf{w})|}{\|\mathbf{w}\|^{3}},
%\end{equation}
%\begin{equation}\kappa(1)= 2 \frac{n-1}{n}\frac{|M^{n-2}D|}{\|\mathbf{w} \|^3}=\frac{n-1}{n}\frac{|\det(M^{n-1}\mathbf{w},M^{n-2}\mathbf{w})|}{\|M^{n-1}\mathbf{w}\|^{3}}
%\end{equation}

\begin{equation}\label{curvature0}
\begin{array}{c}
\kappa(0)%=2 \dfrac{n-1}{n}\dfrac{| D|}{\|\mathbf{w} \|^3}
=\dfrac{n-1}{n}\dfrac{\det(\mathbf{w},M\mathbf{w})}{\|\mathbf{w}\|^{3}},\\[0.5cm]
\kappa(1)%=2 \dfrac{n-1}{n}\dfrac{|M^{n-2}D|}{\|\mathbf{w} \|^3}
=\dfrac{n-1}{n}\dfrac{\det(M^{n-2}\mathbf{w},M^{n-1}\mathbf{w})}{\|M^{n-1}\mathbf{w}\|^{3}}.
\end{array}
\end{equation}

Hence, to calculate the curvature for any $t\in(0,1)$, we reparameterise the arc of the B\'ezier curve defined in $[0,t]$ by subdivision, so this arc is obtained by the action of
the matrix $T=(1-t)\mathbb{I}+tM$ on the vector $t\mathbf{w}$ and the endpoints of this arc are $c(0)$ and $c(t)$. Using the invariance of the curvature
under reparameterisations we get
\begin{equation}\label{curvaturat}
\kappa(t)=\frac{n-1}{n}\frac{\det\left(T^{n-2}\mathbf{w},T^{n-1}\mathbf{w}\right)}{t\|T^{n-1}\mathbf{w} \|^3},
\end{equation}
which gives an expression of the curvature in terms of $\mathbf{w}$, $M$ and $t$.

In order to get a more useful expression of the curvature, it is enough to write the vector $\mathbf{w}$ in a suitable basis on which the action of $M$ on $\mathbf{w}$ is as simple as possible, that is, either a basis of eigenvectors of $M$ or a Jordan basis.

Let $\sigma_1$ and $\sigma_2$ be the eigenvalues of the matrix $M$. They could be real (equal or different) or complex conjugates. If $M$ is a diagonalizable matrix, let $\mathbf{v}_k$ be an eigenvector with eigenvalue $\sigma_k$, $k=1,2$ such that $\{\mathbf{v}_{1},\mathbf{v}_{2}\}$ is a basis of $\mathbb R^2$ (in the case of complex eigenvalues, $\mathbf{v}_2=\overline{\mathbf{v}_1}$). If $M$ does not diagonalize, write $\sigma=\sigma_1=\sigma_2$ and take a Jordan basis $\{\mathbf{v}_{1},\mathbf{v}_{2}\}$ such that $\mathbf{v}_{1}$ is an eigenvector and $\mathbf{v}_{2}$ is a vector satisfying $(M-\sigma\mathbb{I})\mathbf{v}_{2}=\mathbf{v}_{1}$. In this case $\mathbf{v}_1$ and $\mathbf{v}_2$ can be taken to be orthogonal.

Now, any vector $\mathbf{w}\in\mathbb{R}^2$, $\mathbf{w}\neq \mathbf{0}$ can be written with respect to the basis of eigenvectors of $M$ or with respect to the Jordan basis as
\[
\mathbf{w}=\mu_{1}\mathbf{v}_{1}+\mu_{2}\mathbf{v}_{2},
\]
where $\mu_1,\mu_2\in\mathbb R$ if $M$ has real eigenvalues or otherwise,
$\mu_1\in\mathbb{C}\setminus{\mathbb R}$ and $\mu_2=\overline{\mu_1}$.

If $M$ is diagonalizable, then  for $j\ge1$,
\[
M^{j}\mathbf{v}_{1}=\sigma_{1}^{j}\mathbf{v}_{1},\qquad
M^{j}\mathbf{v}_{2}=\sigma_{2}^{j}\mathbf{v}_{2},
\]
and thus,
\[
M^{j}\mathbf{w}=\mu_{1}\sigma_1^{j}\mathbf{v}_{1}+\mu_{2}\sigma_2^{j}\mathbf{v}_{2}.
\]
Written in this way, we get,
\[
\det\left(\mathbf{w}, M\mathbf{w}\right)=\mu_{1}\mu_{2}(\sigma_{2}-\sigma_{1})\det(\mathbf{v}_{1},\mathbf{v}_{2}).
\]
Notice that even in the case of complex eigenvalues $\det\left(\mathbf{w},M\mathbf{w}\right)\in\mathbb{R}$ since $\sigma_{2}-\sigma_{1}=-2i\operatorname{Im}\sigma_1$ and
$\det(\mathbf{v}_{1},\mathbf{v}_{2})=-2i\det\bigl(\operatorname{Re}(\mathbf{v}_{1}),\operatorname{Im}(\mathbf{v}_{1})\bigr)$.

When $M$ does not diagonalize, denote by $\sigma$ the unique eigenvalue of $M$. Then,
\[
M^{j}\mathbf{v}_{1}=\sigma^{j}\mathbf{v}_{1},\qquad
M^{j}\mathbf{v}_{2}=j\sigma^{j-1}\mathbf{v}_{1}+\sigma^{j}\mathbf{v}_{2},
\]
and therefore,
\[
M^{j}\mathbf{w}=\left(\mu_{1}\sigma^{j}+j\mu_{2}\sigma^{j-1}\right)\mathbf{v}_{1}+\mu_{2}\sigma^{j}\mathbf{v}_{2}
\]
which gives in this case,
\[
\det\left(\mathbf{w},M\mathbf{w}\right)=-\mu_{2}^{2}\det(\mathbf{v}_{1},\mathbf{v}_{2}).
\]

In every case, using the properties of the determinant, for $j\ge 1$
\[
\det\left(M^{j-1}\mathbf{w},M^{j}\mathbf{w}\right)=(\sigma_1\sigma_2)^{j-1}\det\left(\mathbf{w},M\mathbf{w}\right),
\]
where $\sigma_1=\sigma_2=\sigma$ if $M$ is not diagonalizable or $\sigma_2=\overline{\sigma_1}$ if the eigenvalues are non-real complex numbers.

Hence by (\ref{curvature0}),
\[
\kappa(1)=\kappa(0)\dfrac{(\sigma_{1}
\sigma_{2})^{n-2}\Vert \mathbf{w}\Vert^3}{\|M^{n-1}\mathbf{w}\|^{3}}.
\]
where
\[
\kappa(0)=\dfrac{n-1}{n}\dfrac{\mu_{1}\mu_{2}(\sigma_{2}-\sigma_{1})\det(\mathbf{v}_{1},\mathbf{v}_{2})}{\|\mathbf{w}\|^{3}}
\]
if $M$ is diagonalizable, or otherwise,
\[
\kappa(0)=-\frac{n-1}{n}\frac{\mu_{2}^{2}\,\det(\mathbf{v}_{1},\mathbf{v}_{2})}
{\|\mathbf{w}\|^{3}}.
\]

Recall that for any $t\in[0,1]$ the curvature is given by the action of $T=(1-t)\mathbb{I}+tM$ on the vector $t\mathbf{w}$. $T$ diagonalizes or not depending on whether $M$ does. If $M$ is diagonalizable then $T$ has eigenvalues $\sigma_{k}(t)=(1-t)+t\sigma_{k}$ (real or complex conjugates) and eigenvectors $\mathbf{v}_k$, $k=1,2$  (that are also the eigenvectors of $M$). Otherwise, the eigenvalues of $T$ are $\sigma_{1}(t)=\sigma_{2}(t)$ where $\{\mathbf{v}_1,\mathbf{v}_2\}$ is the Jordan basis of $M$ chosen above. In every case, by (\ref{curvaturat}),
\[
\kappa(t)=\kappa(0)\dfrac{\bigl(\sigma_1(t)\sigma_2(t)\bigr)^{n-2}\Vert\mathbf{w}\Vert^3}{\Vert T^{n-1}\mathbf{w}\Vert^3},
\]
being $\kappa(0)$ as above.

\

Thus, it has been shown the following
\begin{theorem}\label{tma curvatura}
The curvature of a planar B\'ezier curve of degree $n\ge 2$ gene\-ra\-ted by a matrix $M$ and a vector $\mathbf{w}\neq\mathbf{0}$ is given by
\begin{equation}\label{formula_curvatura}
\kappa(t)=\kappa(0)\dfrac{\bigl(\sigma_1(t)\sigma_2(t)\bigr)^{n-2}\Vert\mathbf{w}\Vert^3}{\Vert T^{n-1}\mathbf{w}\Vert^3},
\end{equation}
where $\kappa(0)$ is the curvature at $t=0$, $T$ is the matrix $T=(1-t)\mathbb{I}+tM$, and for $k=1,2$, $\sigma_{k}(t)=(1-t)+t\sigma_k$ are the eigenvalues of $T$, with $\sigma_k$ the eigenvalues of $M$.

Moreover, writing $\mathbf{w}$ in terms of a basis of eigenvectors of $M$ or a Jordan basis of $M$, as $\mathbf{w}=\mu_1\mathbf{v}_1+\mu_2\mathbf{v}_2$, if $M$ is diagonalizable
then
\[
\kappa(0)=\dfrac{n-1}{n}\dfrac{\mu_{1}\mu_{2}(\sigma_{2}-\sigma_{1})\det(\mathbf{v}_{1},\mathbf{v}_{2})}{\|\mathbf{w}\|^{3}},
\]
otherwise,
\[
\kappa(0)=-\frac{n-1}{n}\frac{\mu_{2}^{2}\,\det(\mathbf{v}_{1},\mathbf{v}_{2})}
{\|\mathbf{w}\|^{3}},
\]
where $\mathbf{v}_{1}$ is an eigenvector of $M$ and $\mathbf{v}_{2}$ satisfies $(M-\sigma\mathbb{I})\mathbf{v}_{2}=\mathbf{v}_{1}$.
\end{theorem}
Notice that the curvature identically vanishes when $\kappa(0)=0$
which occurs when $\mathbf{w}$ is an eigenvector of $M$ and the
resulting B\'ezier curve is a line segment.  On the other hand, if $M$
has complex non-real eigenvalues $\kappa(0)\neq 0$ and hence the
curvature never vanishes regardless the choice of $\mathbf{w}$.
Notice also, that if $n$ is even, or the real eigenvalues of M are
equal or positive, or $M$ has complex eigenvalues, then for any
$t\in(0,1]$, the sign of $\kappa(t)$ is given by the sign of
$\kappa(0)$.

\

Taking the derivative in (\ref{formula_curvatura}) using that $\Vert T^{n-1}\mathbf{w}\Vert^3=\left(\Vert T^{n-1}\mathbf{w}\Vert^2\right)^{3/2}$,
\begin{eqnarray}\label{dcurvature}
\kappa'(t)&=&\dfrac{\kappa(0)\Vert \mathbf{w}\Vert ^{3}\bigl(\sigma_1(t)\sigma_2(t)\bigr)^{n-3}}{2\Vert T^{n-1}\mathbf{w}
\Vert^5}\nonumber\\
&&\cdot
\Bigl(2(n-2)\Vert T^{n-1}\mathbf{w}\Vert^2\bigl(\sigma_1(t)\sigma_2(t)\bigr)'-3\bigl(\sigma_1(t)\sigma_2(t)\bigr)\bigl(\Vert T^{n-1}\mathbf{w}\Vert^2\bigr)'
\Bigr).\nonumber\\
&&
\end{eqnarray}
The expression (\ref{dcurvature}) will be used in the next section to find conditions on $M$ to obtain B\'ezier curves with monotonic curvature.

\section{Conditions for monotonic curvature\label{newconditions}}

We start this section assuming that $M$ has two positive real eigenvalues $\sigma_{1}\ge \sigma_{2}>0$. Consider a basis $\{\mathbf{v}_{1},\mathbf{v}_{2}\}$ of unitary eigenvectors of $M$ such that
$\mathbf{v}_{1}$ is an eigenvector for $\sigma_{1}$ and $\mathbf{v}_{2}$ is an eigenvector for $\sigma_{2}$. In this situation (\ref{dcurvature}) can be rewritten as
\begin{eqnarray}\label{dcurvaturereal}
\kappa'(t)&=&\kappa(0)\dfrac{\Vert \textbf{w}\Vert ^{3}}{2\Vert T^{n-1}\mathbf{w}
\Vert^5}\bigl(\sigma_1(t)\sigma_2(t)\bigr)^{n-3}\nonumber
\\ &\cdot&
\Bigl(-(n+1)\left((\sigma_1+\sigma_2-2)(1-t)+
(2\sigma_1\sigma_2-\sigma_1-\sigma_2)t\right)
\Vert T^{n-1}\mathbf{w}\Vert^2\Bigr.\nonumber\\
&& \Bigl.
-3(n-1)(\sigma_1-\sigma_2)\left(\mu_1^2\sigma_1^{2n-2}(t)-\mu_2^2
\sigma_2^{2n-2}(t)\right)\Bigr)
\end{eqnarray}
where it has been used that
\[
\Vert T^{n-1}\mathbf{w}\Vert^2=\mu_1^2\sigma_1^{2(n-1)}(t)+\mu_2^2\sigma_2^{2(n-1)}(t)+2\mu_1\mu_2\bigl(\sigma_1(t)\sigma_2(t)\bigr)^{n-1}\mathbf{v}_1\cdot\mathbf{v}_2
\]
with $\mathbf{v}_1\cdot\mathbf{v}_2$ the scalar product of $\mathbf{v}_1$ and $\mathbf{v}_2$.

From (\ref{dcurvaturereal}), we give an alternative proof of a result of Cao and Wang (see \cite{cao}) where the matrix $M$ is assumed to be symmetric.

%\begin{corollary} \label{cao}
\begin{theorem}[Cao-Wang] \label{cao}
If the matrix $M$ is symmetric and its eigenvalues,
$\sigma_1\ge\sigma_2>0$, satisfy
\begin{equation}
\label{caowang}
\sigma_{1}\ge1,\qquad 2\sigma_2\ge \sigma_{1}+1,
\end{equation}
then the B\'ezier curve of degree $n\ge 2$ generated by $M$ and $\mathbf{w}$ has monotonic curvature (decreasing if $\kappa(0)>0$, and increasing if $\kappa(0)<0$).
%then the B\'ezier curve of degree $n\ge 3$ generated by $M$ and $\mathbf{w}$ has decreasing curvature.
%\end{corollary}
\end{theorem}

\

\noindent {\it Proof}. Since $\sigma_1\ge \sigma_2$ conditions in (\ref{caowang}) imply also that
\begin{equation}\label{caowang2}
\sigma_2\ge 1,\quad \text{ and }\quad 2\sigma_1\ge \sigma_2+1.
\end{equation}
When $M$ is symmetric its eigenvectors $\mathbf{v}_1$ and $\mathbf{v}_2$ are orthogonal, and then,
\[
\Vert T^{n-1}\mathbf{w}\Vert^2=\mu_1^2\sigma_1^{2n-2}(t)+\mu_2^2\sigma_2^{2n-2}(t).
\]
It is easy to check that $k'(t)$ does not change its sign for $t\in[0,1]$. Indeed, we may group the large parenthesis of the derivative of the
curvature (\ref{dcurvaturereal}) so the terms that multiply $\mu_1^2\sigma_1^{2n-2}(t)$ and $\mu_2^2\sigma_2^{2n-2}(t)$ are respectively,
\[
-\Bigl(2(n+1)(\sigma_{1}-1)(\sigma_{2}-1)t+2(2\sigma_{1}-\sigma_{2}-1)n+2(2\sigma_{2}-\sigma_{1}-1) \Bigr),
\]
and
\[
-\Bigl(2(n+1)(\sigma_{1}-1)(\sigma_{2}-1)t+2(2\sigma_{2}-\sigma_{1}-1)n
+2(2\sigma_{1}-\sigma_{2}-1) \Bigr).
\]
Thus the sign of every summand is given by (\ref{caowang}) or (\ref{caowang2})
and hence, $\kappa'(t)/\kappa(0)<0$ for every $t\in[0,1]$ and for any degree $n$ as shown in \cite[Condition (6)]{cao}.
\hfill $\square$

\

For a non-symmetric matrix $M$ we get the following result.

\begin{theorem}\label{ejemplos}
If the eigenvalues of $M$, $\sigma_1\ge\sigma_2>0$, satisfy
\begin{equation} \label{autovalores}
\sigma_1+\sigma_2-2 \ge 0,
\end{equation}
and the vector $\mathbf{w}=\mu_1\mathbf{v}_1+\mu_2\mathbf{v}_2\neq\mathbf{0}$ is chosen so that
\begin{equation} \label{restriccion_vector}
\vert \mu_1 \vert \ge \vert\mu_2\vert > 0,
\end{equation}
then the B\'ezier curve of degree $n\ge 2$ generated by $M$ and $\mathbf{w}$ has monotonic curvature (decreasing if $\kappa(0)>0$ and increasing if $\kappa(0)<0$).
\end{theorem}

When $\sigma_1,\sigma_2\ge 1$ condition (\ref{autovalores}) is
immediately satisfied, so (\ref{autovalores}) only has content when no
a priori lower bound on the eigenvalues (other than their positivity)
is assumed.  On the other hand, condition (\ref{autovalores}) implies
that the largest eigenvalue, $\sigma_1$ is at least $1$, and gives a
lower bound for the smallest eigenvalue, $\sigma_2$, in terms of
$\sigma_1$, namely $\sigma_2\ge 2-\sigma_1$.

If $\sigma_1,\sigma_2\ge 1$, our approach differs from that of Cao and Wang (see \cite{cao}) and of Farin (see \cite{farin-a}), in that the restriction is not applied on the eigenvalues, but rather on the vectors on which the matrix $M$ acts, that is, condition (\ref{restriccion_vector}). Notice that, by the affine invariance property of B\'ezier curves and the fact that the curvature does not change its monotonicity by rotation or scaling, (\ref{restriccion_vector}) does not imply any restriction on the direction of the tangent to the B\'ezier curve at its initial point.

One interesting feature of condition (\ref{autovalores}) is that it is
not invariant under subdivision, so Farin's approach in \cite{farin-a}
could not lead to such a result.  Indeed, the matrix that generates
the second subarc of the B\'ezier curve from $t$ to $1$, $MT^{-1}$
(where $T=(1-t)\mathbb{I}+tM$), has eigenvalues
$\sigma_k/\sigma_k(t)$, $k=1,2$, for which condition
(\ref{autovalores}) does not necessarily hold (for example if
$\sigma_1=2$, $\sigma_2=1/2$ and $t=3/4$).

\

\noindent {\it Proof}. Without loss of generality assume $\sigma_1>\sigma_2$, (that is, the matrix $M$ is not
a constant multiple of the identity) and that the eigenvectors $\mathbf{v}_1$ and $\mathbf{v}_2$ are unitary. The sign of the derivative (\ref{dcurvaturereal}) of the curvature can be obtained from the sign of the term inside the big parenthesis that will be denoted by $P_1$,
\[
\begin{split}
P_1:=&-(n+1)\bigl((\sigma_1+\sigma_2-2)(1-t)+(2\sigma_1\sigma_2-\sigma_1-\sigma_2)t\bigr)\Vert T^{n-1}\mathbf{w}\Vert^2\\
&-3(n-1)(\sigma_1-\sigma_2)\bigl(\mu_1^2\sigma_1^{2n-2}(t)-\mu_2^2\sigma_2^{2n-2}(t)\bigr).
\end{split}
\]
By condition (\ref{restriccion_vector}) the second summand in $P_1$
is always negative.  If for a given $t$, the first summand in $P_1$ is non-positive then $\kappa'(t)/\kappa(0)<0$. Thus we are left to prove that $\kappa'(t)/\kappa(0)<0$ still holds for those $t$ for which the first summand in $P_1$ is positive, that is, $t$ such that $(\sigma_1+\sigma_2-2)(1-t)+(2\sigma_1\sigma_2-\sigma_1-\sigma_2)t<0$. Then, by condition (\ref{autovalores}), necessarily $t>0$ and $\sigma_2< 1$.

So assume that $(2-\sigma_1-\sigma_2)(1-t)+(\sigma_1+\sigma_2-2\sigma_1\sigma_2)t>0$. Since $\mathbf{v}_1$ and $\mathbf{v}_2$ are chosen to be unitary,
\[
\begin{split}
\Vert T^{n-1}\mathbf{w}\Vert^2&=\mu_1^2\sigma_1^{2n-2}(t)+\mu_2^2\sigma_2^{2n-2}(t)+2\mu_1\mu_2\bigl(\sigma_1(t)\sigma_2(t)\bigr)^{n-1}\mathbf{v}_1\cdot\mathbf{v}_2\\
& < \mu_1^2\sigma_1^{2n-2}(t)+\mu_2^2\sigma_2^{2n-2}(t)+2\vert\mu_1\mu_2\vert\bigl(\sigma_1(t)\sigma_2(t)\bigr)^{n-1}\\
&=\left(\vert\mu_1\vert\sigma_1^{n-1}(t)+\vert\mu_2\vert\sigma_2^{n-1}(t)\right)^2.
\end{split}
\]
%and $P_1$ can be bounded as
%\[
%\begin{split}
%P_1&\le %(n+1)\bigl((2-\sigma_1-\sigma_2)(1-t)+(\sigma_1+\sigma_2-2\sigma_1\sigma_2)t\bigr%)\left(\vert\mu_1\vert\sigma_1^{n-1}(t)+\vert\mu_2\vert\sigma_2^{n-1}(t)\right)^2%\\
%&\phantom{=}-3(n-1)(\sigma_1-\sigma_2)\bigl(\mu_1^2\sigma_1^{2n-2}(t)-\mu_2^2\sig%ma_2^{2n-2}(t)\bigr).
%\end{split}
%\]
Using this bound in the expression of $P_1$, taking common factor 
\[
\vert\mu_1\vert\sigma_1^{n-1}(t)+\vert\mu_2\vert\sigma_2^{n-1}(t)
\] 
and recalling that $\vert\mu_1\vert \geq \vert\mu_2\vert$ (condition (\ref{restriccion_vector})),
\[
\begin{split}
P_1&\le  \vert\mu_1\vert\left(\vert\mu_1\vert\sigma_1^{n-1}(t)+\vert\mu_2\vert\sigma_2^{n-1}(t)\right)\cdot\\
&\phantom{=} \Bigl((n+1)\bigl((2-\sigma_1-\sigma_2)(1-t)+(\sigma_1+\sigma_2-2\sigma_1\sigma_2)t\bigr)\left(\sigma_1^{n-1}(t)+\sigma_2^{n-1}(t)\right)\\
&\phantom{=(}-3(n-1)(\sigma_1-\sigma_2)\bigl(\sigma_1^{n-1}(t)-\sigma_2^{n-1}(t)\bigr)\Bigr),
\end{split}
\]
so again, the sign of $P_1$ is given by the sign of the factor inside the big parenthesis,
\[
\begin{split}
P_2:=&(n+1)\bigl((2-\sigma_1-\sigma_2)(1-t)+(\sigma_1+\sigma_2-2\sigma_1\sigma_2)t\bigr)\left(\sigma_1^{n-1}(t)+\sigma_2^{n-1}(t)\right)\\
&-3(n-1)(\sigma_1-\sigma_2)\bigl(\sigma_1^{n-1}(t)-\sigma_2^{n-1}(t)\bigr).
\end{split}
\]
Consider $P_2$ as a function of $n$, that is $P_2=P_2(n)$. Since $\sigma_1>\sigma_2>0$ then $\sigma_1(t)> \sigma_2(t)>0$, and recall that we are under the assumption that
$(2-\sigma_1-\sigma_2)(1-t)+(\sigma_1+\sigma_2-2\sigma_1\sigma_2)t>0$. Therefore, denoting 
\[A=(n+1)\bigl((2-\sigma_1-\sigma_2)(1-t)+(\sigma_1+\sigma_2-2\sigma_1\sigma_2)t\bigr)\left(\sigma_1^{n-1}(t)+\sigma_2^{n-1}(t)\right),\] 
\[B=-3(n-1)(\sigma_1-\sigma_2)\bigl(\sigma_1^{n-1}(t)-\sigma_2^{n-1}(t)\bigr),\]
both $A$ and $B$ are positive, $P_2(n)=A-B$ and we can consider the
function $f(n)$ given by $f(n)=\log A - \log B$. 

Clearly, $P_2(n)<0$ if and only if $f(n)<0$. Since $f(n)$ is a decreasing function of $n$ for $n\ge 2$, it is enough to show that $P_2(2)<0$ to obtain that $P_2=P_2(n)<0$ for every $n\ge 2$.
Thus, taking $n=2$ 
\[
\begin{split}
P_2(2)= &3\bigl((2-\sigma_1-\sigma_2)(1-t)+(\sigma_1+\sigma_2-2\sigma_1\sigma_2)t\bigr)\left(\sigma_1(t)+\sigma_2(t)\right)\\
&-3(\sigma_1-\sigma_2)\bigl(\sigma_1(t)-\sigma_2(t)\bigr).
\end{split}
\]
Substituting $\sigma_1(t)+\sigma_2(t)=2(1-t)+t(\sigma_1+\sigma_2)$ and $\sigma_1(t)-\sigma_2(t)=t(\sigma_1-\sigma_2)$, the expression above can be written as
\[
\begin{split}
P_2(2)= & \, 3(2-\sigma_1-\sigma_2)(1-t)\left(2(1-t)+t(\sigma_1+\sigma_2)\right)\\
&+3t(\sigma_1+\sigma_2-2\sigma_1\sigma_2)\left(2(1-t)+t(\sigma_1+\sigma_2)\right)\\
&-3t(\sigma_1-\sigma_2)^2.
\end{split}
\]
Combining the terms with common factor $3t$, writing $(\sigma_1-\sigma_2)^2=(1-t)(\sigma_1-\sigma_2)^2+t(\sigma_1-\sigma_2)^2$ and reorganizing the terms on the right hand side of the inequality we get
\[
\begin{split}
P_2(2)= & \, 3(2-\sigma_1-\sigma_2)(1-t)\left(2(1-t)+t(\sigma_1+\sigma_2)\right)\\
&+3t(1-t)(\sigma_1+\sigma_2)(2-\sigma_1-\sigma_2)+6t^2\sigma_1\sigma_2(2-\sigma_1-\sigma_2).
\end{split}
\]
By (\ref{autovalores}), the common factor $2-\sigma_1-\sigma_2$ is negative, and every other factor is positive therefore $P_2(2)<0$ which gives that $P_2<0$ for any $n\ge 2$ which in turn implies that $\kappa'(t)/\kappa(0)<0$ for any $n\ge 2$, even for those $t$ such that $(\sigma_1+\sigma_2-2)(1-t)+(2\sigma_1\sigma_2-\sigma_1-\sigma_2)t<0$.

Hence $\kappa'(t)/\kappa(0)<0$ for every $t\in[0,1]$, and therefore $\kappa(t)$ is monotonic in $[0,1]$. \hfill $\square$

\

\noindent {\bf Remark}.  Although condition (\ref{autovalores}) is not
invariant by subdivision the stronger condition
$2\sigma_1\sigma_2-\sigma_1-\sigma_2\ge 0$ that appears in the
derivative of the curvature is.

Below, we include some examples to show that without conditions
(\ref{autovalores}) and (\ref{restriccion_vector}) in Theorem
\ref{ejemplos}, the resulting curve may not have monotonic curvature.
The first example gives a counterexample when condition
(\ref{autovalores}) does not hold.

\begin{example}
Let $M_{1}=\left(\begin{smallmatrix}5/4 & 0 \\ 0 & 1/10\end{smallmatrix}\right)$
and let $\mathbf{w}_1=\left(\begin{smallmatrix}1\\-1\end{smallmatrix}\right)$. Then the cubic B\'ezier curve generated by $M_1$ and $\mathbf{w}_1$ does not have monotonic curvature. See Figure \ref{eje1}.

\begin{figure}[h]
\centering
\begin{subfigure}{0.4\textwidth}
\includegraphics[width=\textwidth]{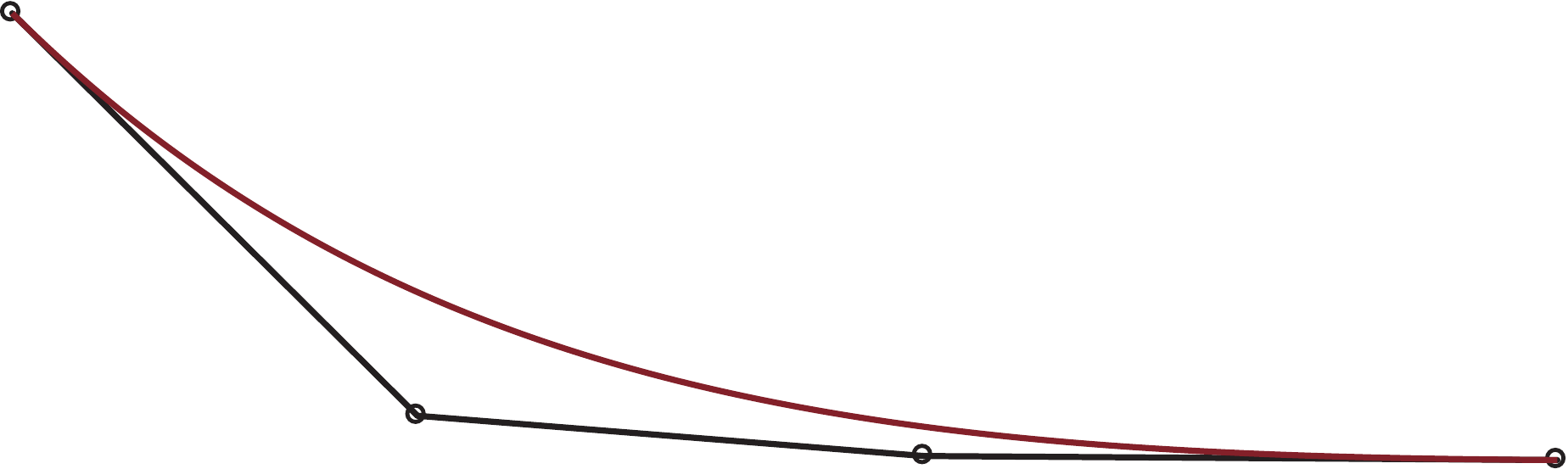}
\end{subfigure}
\hspace{1cm}
\begin{subfigure}{0.3\textwidth}
\includegraphics[width=\textwidth]{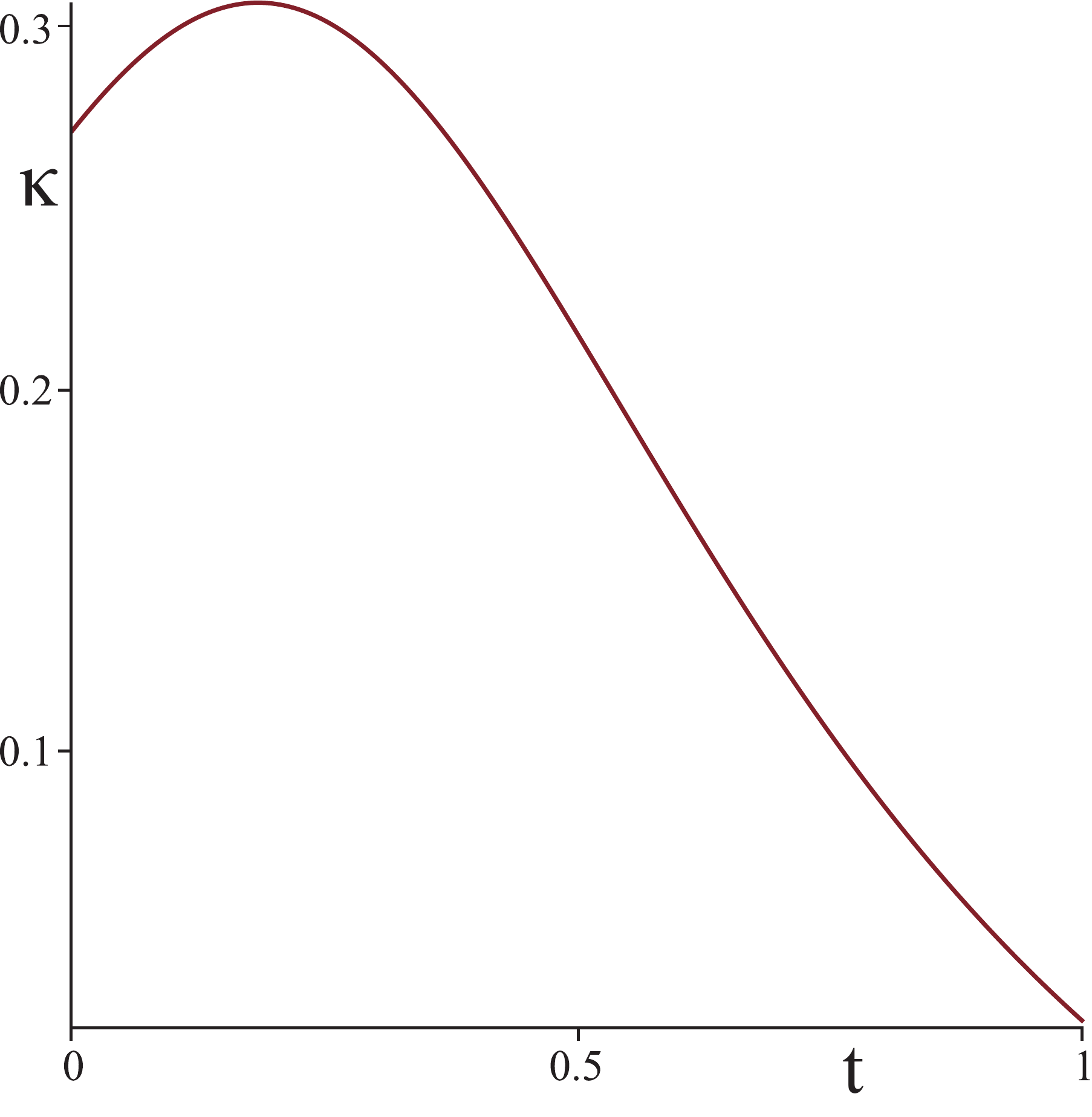}
\end{subfigure}
\caption{Cubic B\'ezier curve generated by $M_1$ and $\mathbf{w}_1$ and the graph of its curvature}\label{eje1}
\end{figure}
\end{example}

The next example is a counterexample of monotonic curvature when condition
(\ref{restriccion_vector}) is not satisfied

\begin{example}
Let $M_2$ be the matrix $M_2=\left(\begin{smallmatrix}4 & 0 \\ 0 & 1\end{smallmatrix}\right)$
and let $\mathbf{w}_2=\left(\begin{smallmatrix}2/5 \\- 5\end{smallmatrix}\right)$. Then the cubic B\'ezier curve generated by $M_2$ and $\mathbf{w}_2$ does not have monotonic curvature. See Figure \ref{eje2}.

\begin{figure}[h]
\centering
\begin{subfigure}{0.4\textwidth}
\includegraphics[height=4cm]{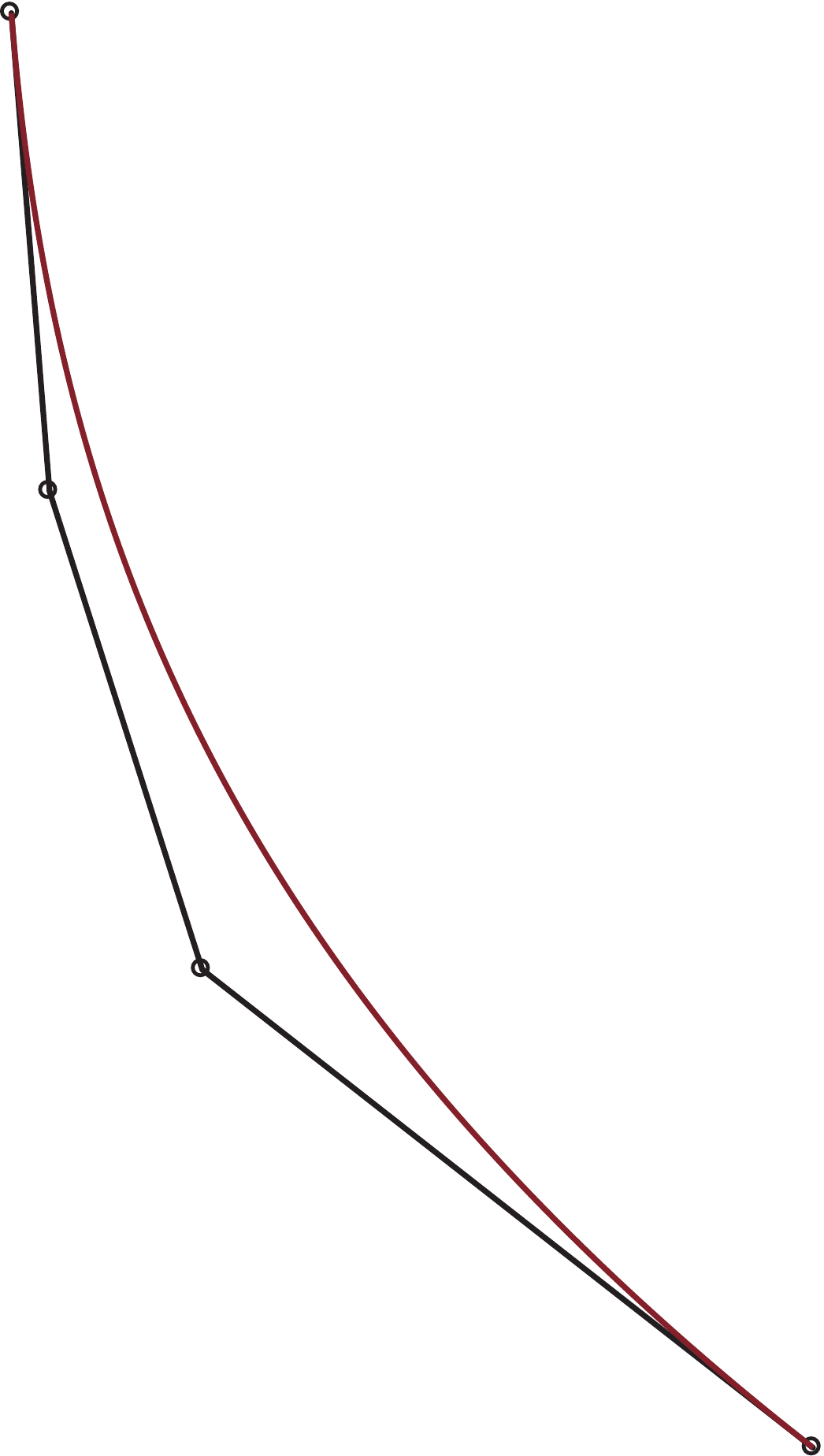}
\end{subfigure}
\hspace{1cm}
\begin{subfigure}{0.3\textwidth}
\includegraphics[width=\textwidth]{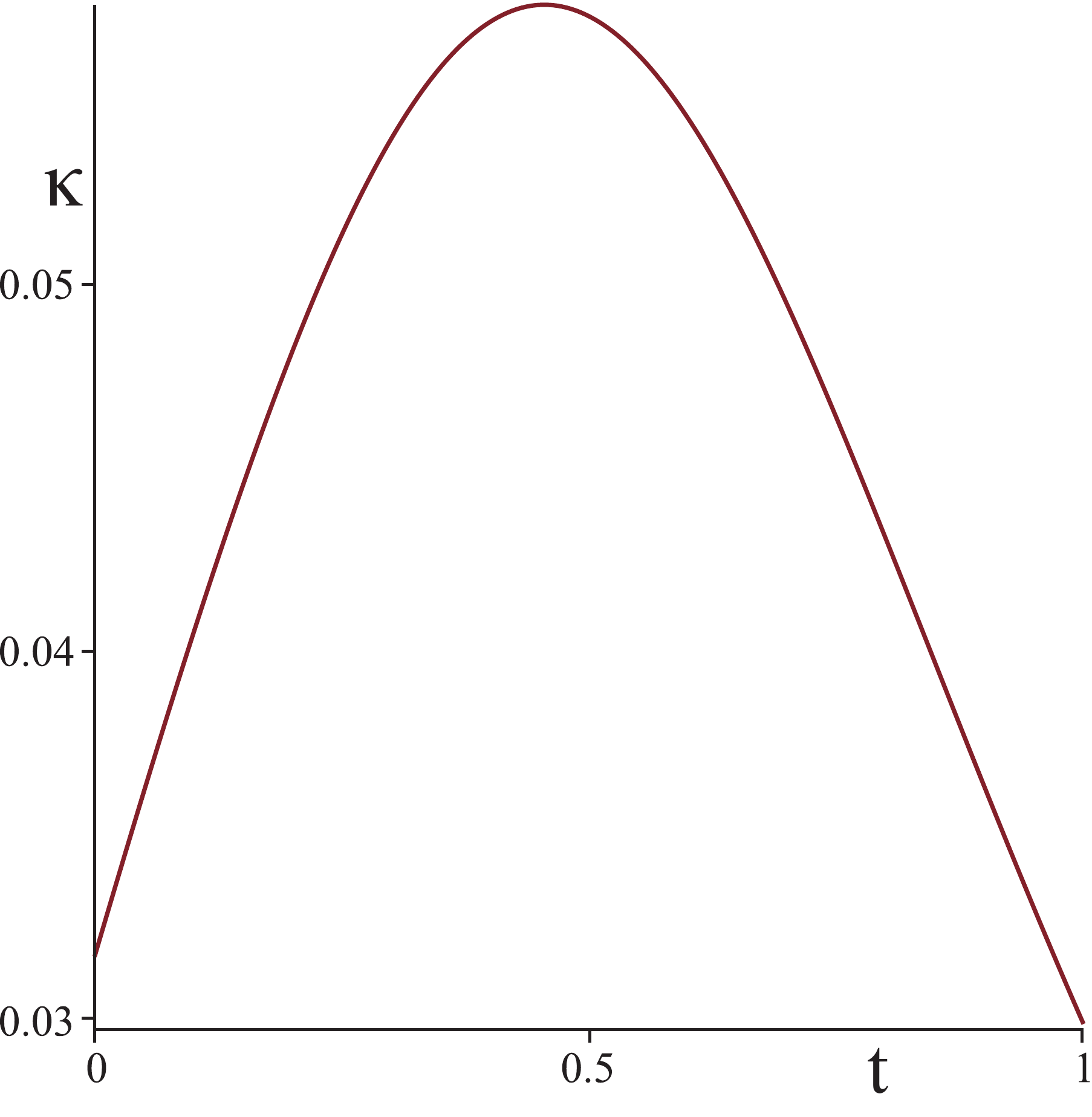}
\end{subfigure}
\caption{Cubic B\'ezier curve generated by $M_2$ and $\mathbf{w}_2$ and the graph of its curvature}\label{eje2}
\end{figure}
\end{example}

And finally we show a quintic B\'ezier curve generated by matrix and vector as in Theorem \ref{ejemplos}.

\begin{example}
Consider the vectors
$\mathbf{v}_1=\left(\begin{smallmatrix}1\\0\end{smallmatrix}\right)$
and $\mathbf{v}_2=\left(\begin{smallmatrix}-\sqrt{3}/2 \\
-1/2\end{smallmatrix}\right)$ so that the angle $\gamma$ between
$\mathbf{v}_1$ and $\mathbf{v}_2$ is $\gamma=-5\pi/6$.  Let $M_3$ be
the matrix with eigenvectors $\mathbf{v}_1$ and $\mathbf{v}_2$ with
eigenvalues $\sigma_1=3/2$ and $\sigma_2=3/4$ (so condition
(\ref{autovalores}) holds for these eigenvalues).  Take
$\mathbf{w}_3=2\mathbf{v}_1+2\mathbf{v}_2$.  Hence
$M_{3}=\left(\begin{smallmatrix}3/2 & -3\sqrt{3}/4 \\ 0 &
3/4\end{smallmatrix}\right)$ and
$\mathbf{w}_3=\left(\begin{smallmatrix}2-\sqrt{3}\\-1\end{smallmatrix}\right)$.

By Theorem \ref{ejemplos} the B\'ezier curve of degree 5 generated by $M_3$ and $\mathbf{w}_3$ has monotonic curvature. See Figure \ref{eje3}.

\begin{figure}[h]
\centering
\begin{subfigure}{0.55\textwidth}
\includegraphics[width=\textwidth]{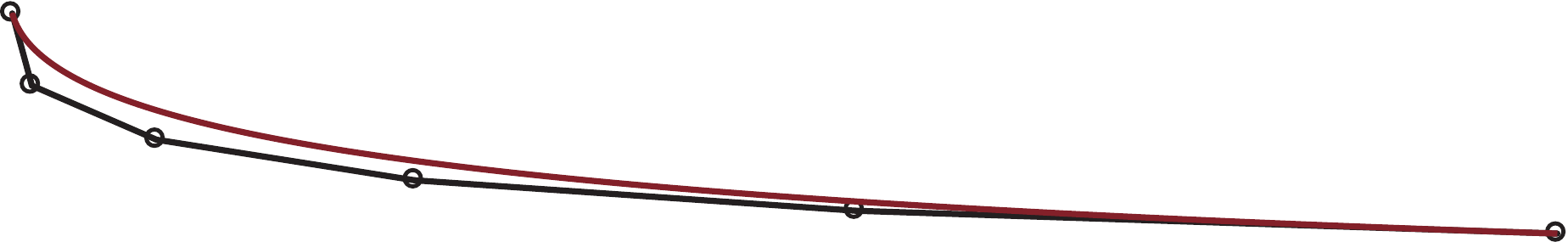}
\end{subfigure}
\begin{subfigure}{0.35\textwidth}
\includegraphics[width=\textwidth]{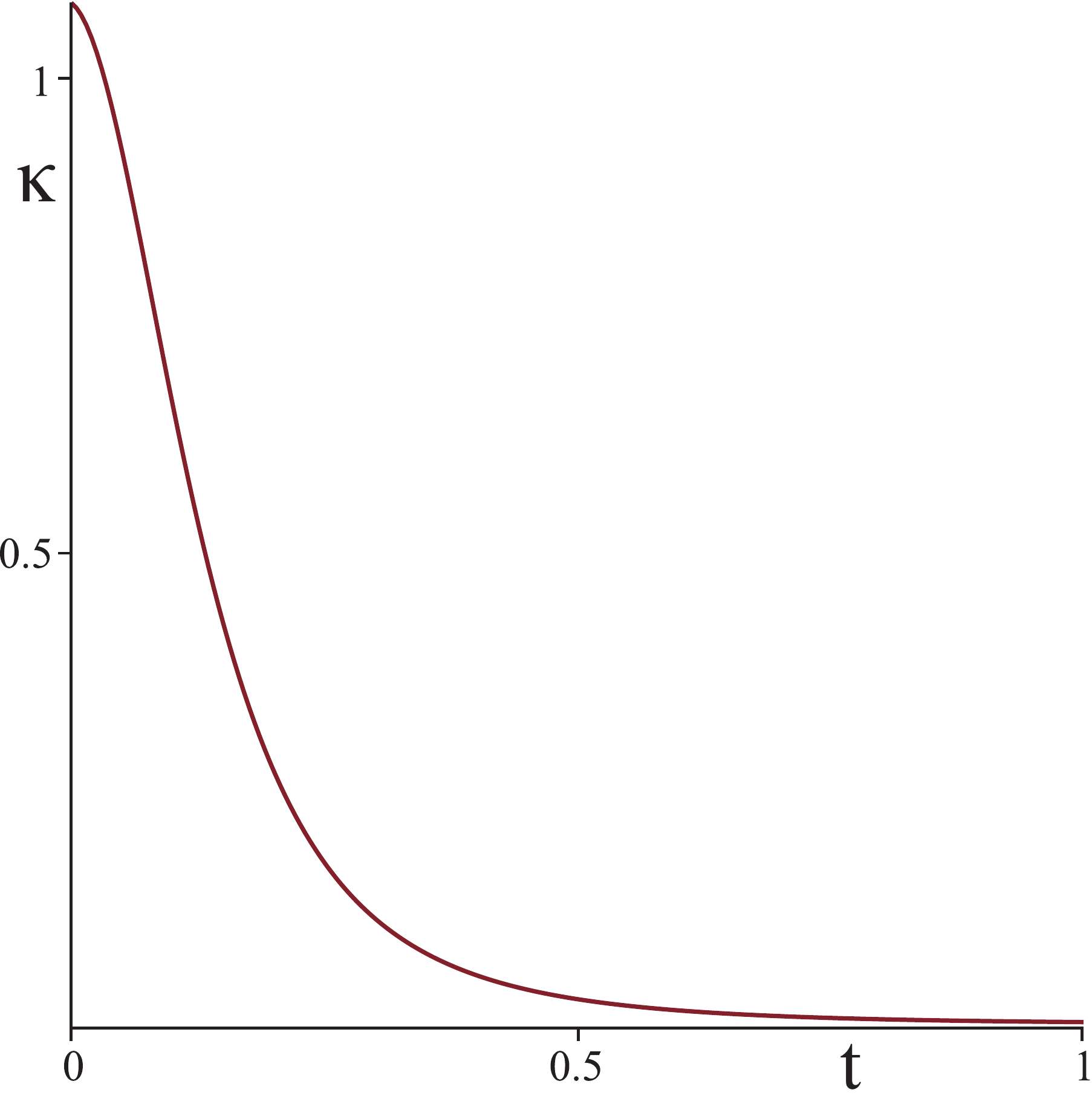}
\end{subfigure}
\caption{Quintic B\'ezier curve generated by $M_3$ and $\mathbf{w}_3$ and the graph of its curvature}\label{eje3}
\end{figure}

\end{example}

%%%% No cumplen 

Conditions (\ref{autovalores}) and (\ref{restriccion_vector}) in Theorem
\ref{ejemplos} are not necessary conditions to obtain curves with 
monotononic curvature as the next two examples show:

\begin{example}
Consider the matrix $M_{4}=\left(\begin{smallmatrix} 3/2 &
6/5\sqrt{3}\\0 & 3/10\end{smallmatrix}\right)$ with
eigenvalues $\sigma_1=3/2$ and $\sigma_2=3/10$, and
corresponding eigenvectors $\mathbf{v}_1=\left(\begin{smallmatrix} 1
\\ 0\end{smallmatrix}\right)$ and
$\mathbf{v}_2=\left(\begin{smallmatrix} 1/2 \\
-\sqrt{3}/2\end{smallmatrix}\right)$. Then the B\'ezier curve of
degree $4$ generated by $M$ and
$\mathbf{w}_{4}=\mathbf{v}_1+\tfrac{1}{2}\mathbf{v}_2$ has monotonic
curvature in spite that condition (\ref{autovalores}) does not hold. 
See Figure \ref{eje12}.
\begin{figure}[h]
\centering
\begin{subfigure}{0.55\textwidth}
\includegraphics[width=\textwidth]{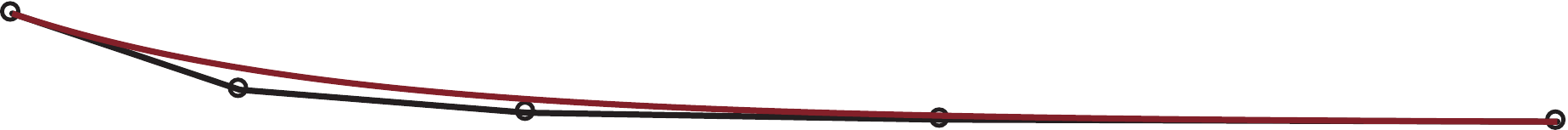}
\end{subfigure}
\begin{subfigure}{0.35\textwidth}
\includegraphics[width=\textwidth]{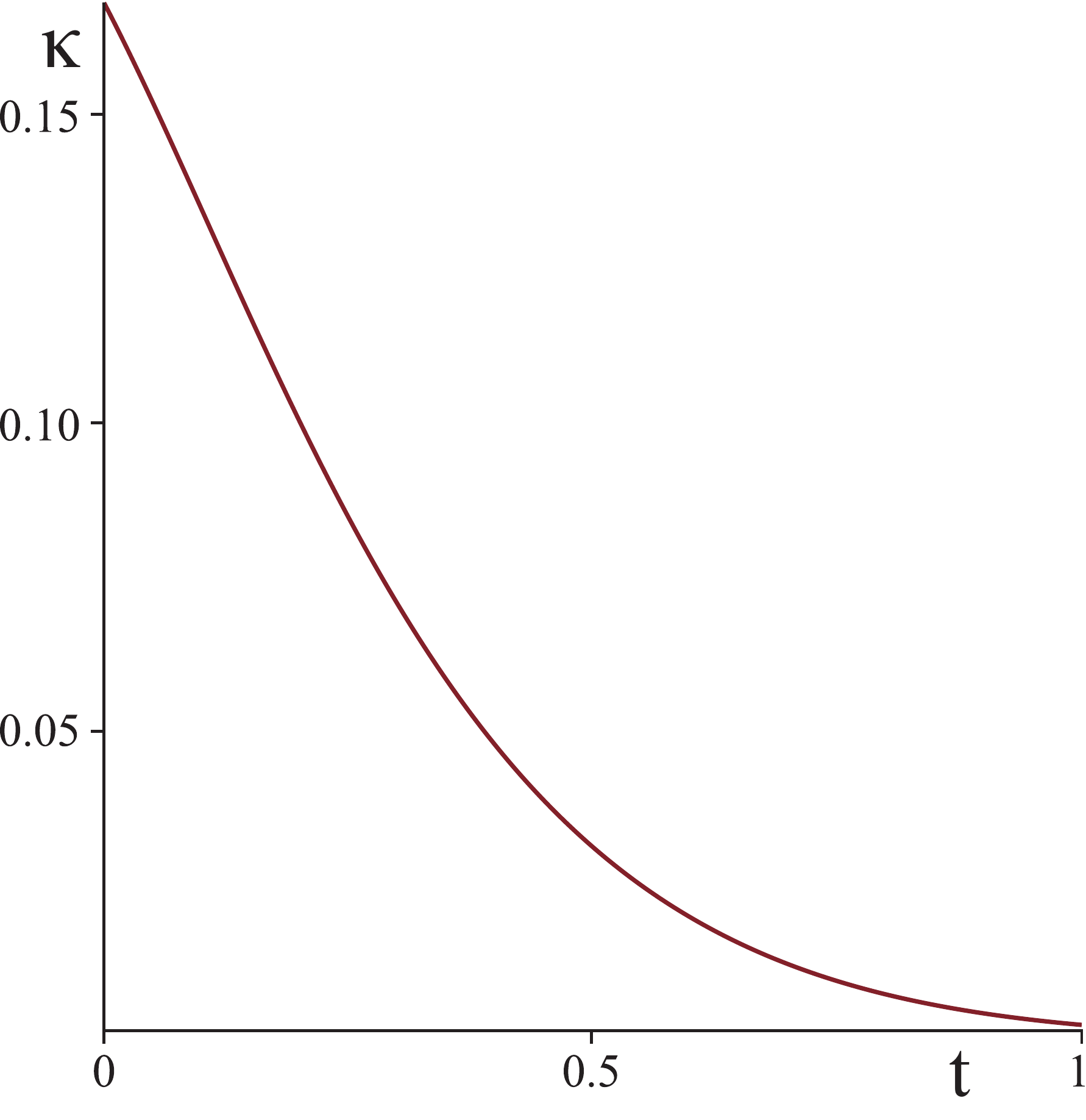}
\end{subfigure}
\caption{Quartic B\'ezier curve generated by $M_4$ and $\mathbf{w}_4$
and the graph of its curvature}\label{eje12}
\end{figure}

\end{example}

\begin{example}
Consider the matrix $M_{5}=\left(\begin{smallmatrix} 3/2 & 0
\\0 & 7/10\end{smallmatrix}\right)$ with eigenvalues
$\sigma_1=3/2$ and $\sigma_2=7/10$, and
corresponding eigenvectors $\mathbf{v}_1=\left(\begin{smallmatrix} 1
\\ 0\end{smallmatrix}\right)$ and
$\mathbf{v}_2=\left(\begin{smallmatrix} 0 \\
1\end{smallmatrix}\right)$. Then the cubic curve generated by $M$ and
$\mathbf{w}_{5}=\mathbf{v}_1-\tfrac{11}{10}\mathbf{v}_2$ has monotonic
curvature in spite that condition (\ref{restriccion_vector}) does not
hold. Notice that this example does not satisfy condition
(\ref{caowang}) in Theorem \ref{cao} of Cao and Wang. See Figure 
\ref{eje13}.
\begin{figure}[h]
\centering
\begin{subfigure}{0.45\textwidth}
\includegraphics[width=\textwidth]{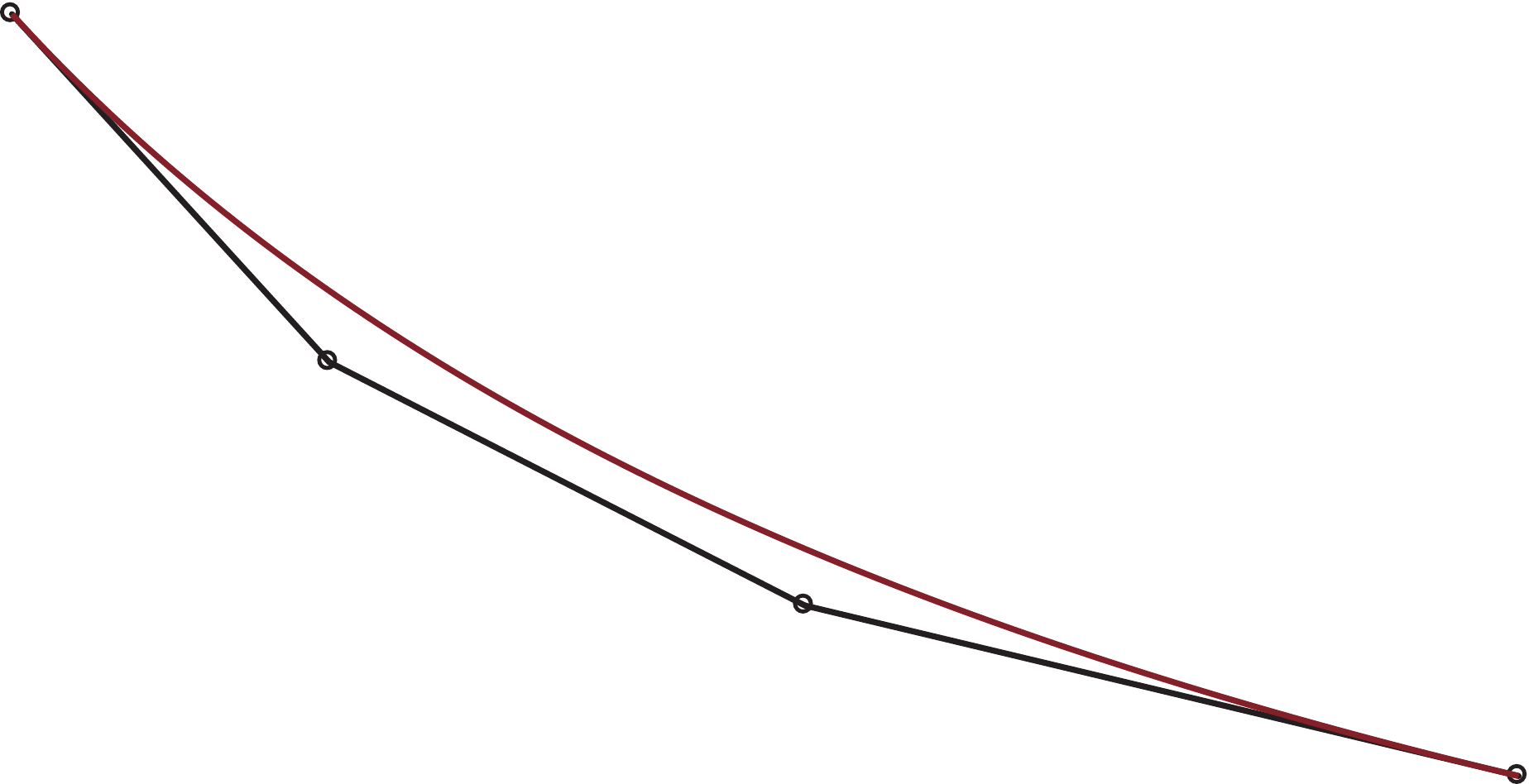}
\end{subfigure}
\begin{subfigure}{0.35\textwidth}
\includegraphics[width=\textwidth]{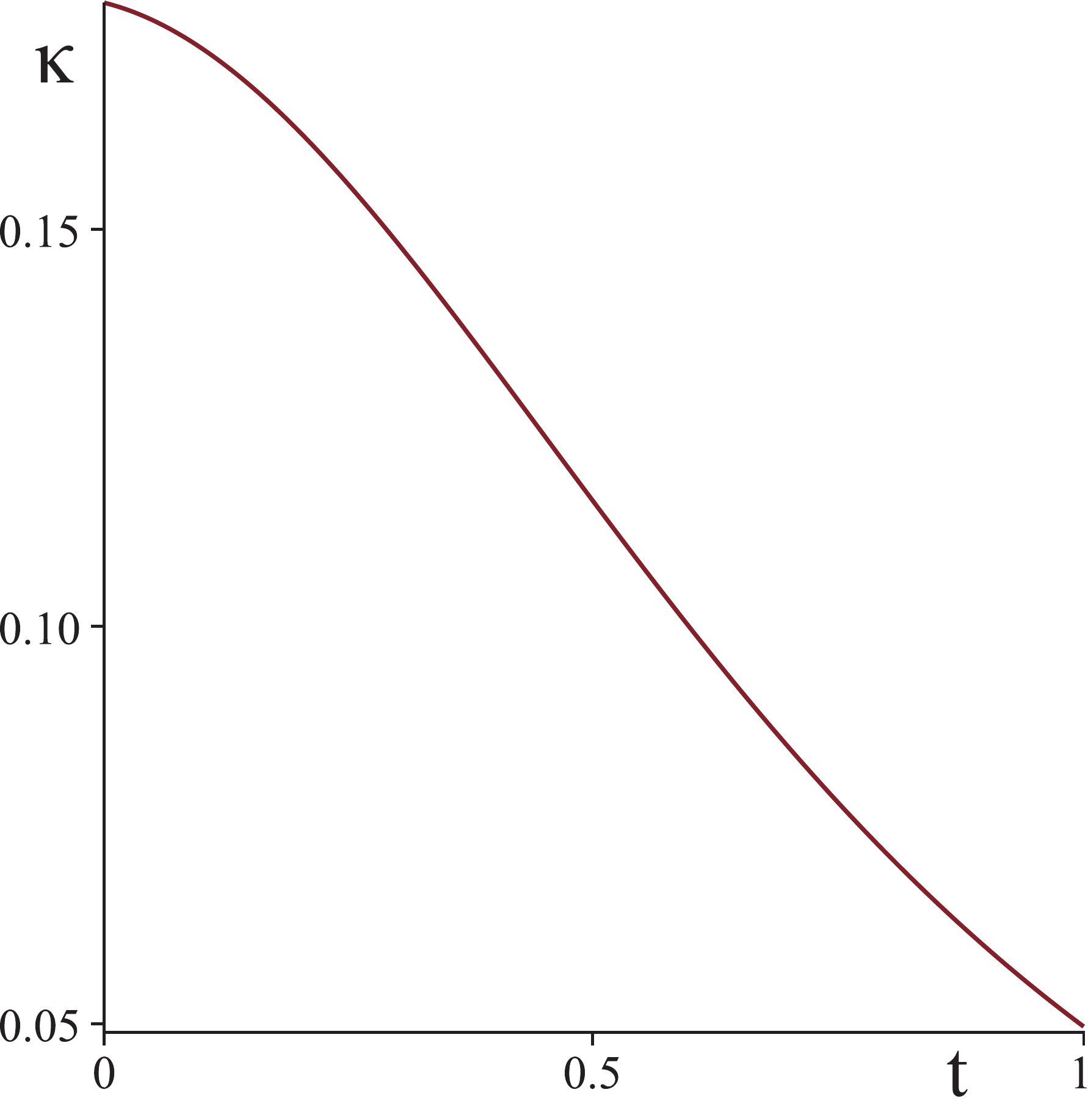}
\end{subfigure}
\caption{Cubic B\'ezier curve generated by $M_5$ and $\mathbf{w}_5$
and the graph of its curvature}\label{eje13}
\end{figure}

\end{example}

Next, we consider the case of a non-diagonalizable matrix $M$ with (double) real eigenvalue $\sigma$. As mentioned above, the Jordan basis $\{\mathbf{v}_1,\mathbf{v}_2\}$ is chosen to be orthogonal where $\mathbf{v}_1$ is an eigenvector of $M$ and $\mathbf{v}_2$ satisfies $(M-\sigma\mathbb{I})\mathbf{v}_{2}=\mathbf{v}_{1}$.

\begin{theorem}\label{tma:Jordan}
If $\sigma$, the eigenvalue of $M$, satisfies $\sigma\ge 1$, and the vector $\mathbf{w}=\mu_1\mathbf{v}_1+\mu_2\mathbf{v}_2\neq\mathbf{0}$ is chosen with coefficients such that
\[
\mu_1\mu_2\ge 0, \quad \mu_2\neq 0,
\]
then the B\'ezier curve of degree $n\ge 2$ generated by $M$ and $\mathbf{w}$ has monotonic curvature (decreasing if $\kappa(0)>0$ and increasing if $\kappa(0)<0$).
\end{theorem}

\

\noindent {\it Proof}. Recall that for $t\in[0,1]$ the matrix $T$ is given by $T=(1-t)\mathbb{I}+tM$ with eigenvalue
\[
\sigma(t)=1-t+t\sigma.
\]
Substituting in (\ref{dcurvature}), $\sigma_1(t)=\sigma_2(t)=\sigma(t)$ and,
\[
\Vert T^{n-1}\mathbf{w}\Vert^2=\bigl(\mu_1\sigma^{n-1}(t)+\mu_2(n-1)t\sigma^{n-2}(t)\bigr)^2\Vert\mathbf{v}_1\Vert^2+\bigl(\mu_2\sigma^{n-1}(t)\bigr)^2\Vert\mathbf{v}_2\Vert^2,
\]
(where the orthogonality of $\mathbf{v}_1$ and $\mathbf{v}_2$ has been used), we get
\[
\begin{split}
\kappa'(t)=-\dfrac{\kappa(0)\Vert \mathbf{w}\Vert^3}{\sigma^{n-1}(t)}&\left(\dfrac{(n+1)(\sigma-1)}{\bigl((\mu_1\sigma(t)+\mu_2(n-1)t)^2\Vert\mathbf{v}_1\Vert^2+(\mu_2\sigma(t))^2\Vert\mathbf{v}_2\Vert^2\bigr)^{3/2}}\right.\\[0.2cm]
&\hspace{0.2cm}\left.+\dfrac{3\mu_2(n-1)\bigl(\mu_1\sigma(t)+\mu_2(n-1)t\bigr)\Vert\mathbf{v}_1\Vert^2}
{\bigl((\mu_1\sigma(t)+\mu_2(n-1)t)^2\Vert\mathbf{v}_1\Vert^2+(\mu_2\sigma(t))^2\Vert\mathbf{v}_2\Vert^2\bigr)^{5/2}}\right).
\end{split}
\]

Since $\sigma\ge 1$, $\mu_1\mu_2\ge 0$ and $\mu_2\neq 0$, every term in the parenthesis is positive for any $t\in[0,1]$, and hence $\kappa'(t)/\kappa(0)<0$ in $[0,1]$.  \hfill $\square$

\

In the following examples, the Jordan basis
$\{\mathbf{v}_1,\mathbf{v}_2\}$ is given by $\mathbf{v}_1 =
\left(\begin{smallmatrix}0\\1\end{smallmatrix} \right)$ and
$\mathbf{v}_2 = \left(\begin{smallmatrix}1\\0\end{smallmatrix} \right)$.  The
next two examples show that without conditions $\sigma\ge 1$ or
$\mu_{1}\mu_{2}\ge 0$ the curve may not have monotonic curvature.

% The next examples show the
% necessity of conditions $\sigma \ge 1$ and $\mu_1\mu_2\ge 0$ in the
% case of a non-diagonalizable matrix.

\begin{example}
Let $M_{6}=\left(\begin{smallmatrix}0.7 & 0 \\ 1 & 0.7\end{smallmatrix}\right)$
and $\mathbf{w}_6=\left(\begin{smallmatrix}10\\1\end{smallmatrix}\right)$. 
Then the cubic B\'ezier curve generated by $M_6$ and $\mathbf{w}_6$ does not have monotonic curvature. See Figure \ref{eje4}.

\begin{figure}[h]
\centering
\begin{subfigure}{0.25\textwidth}
\includegraphics[width=\textwidth]{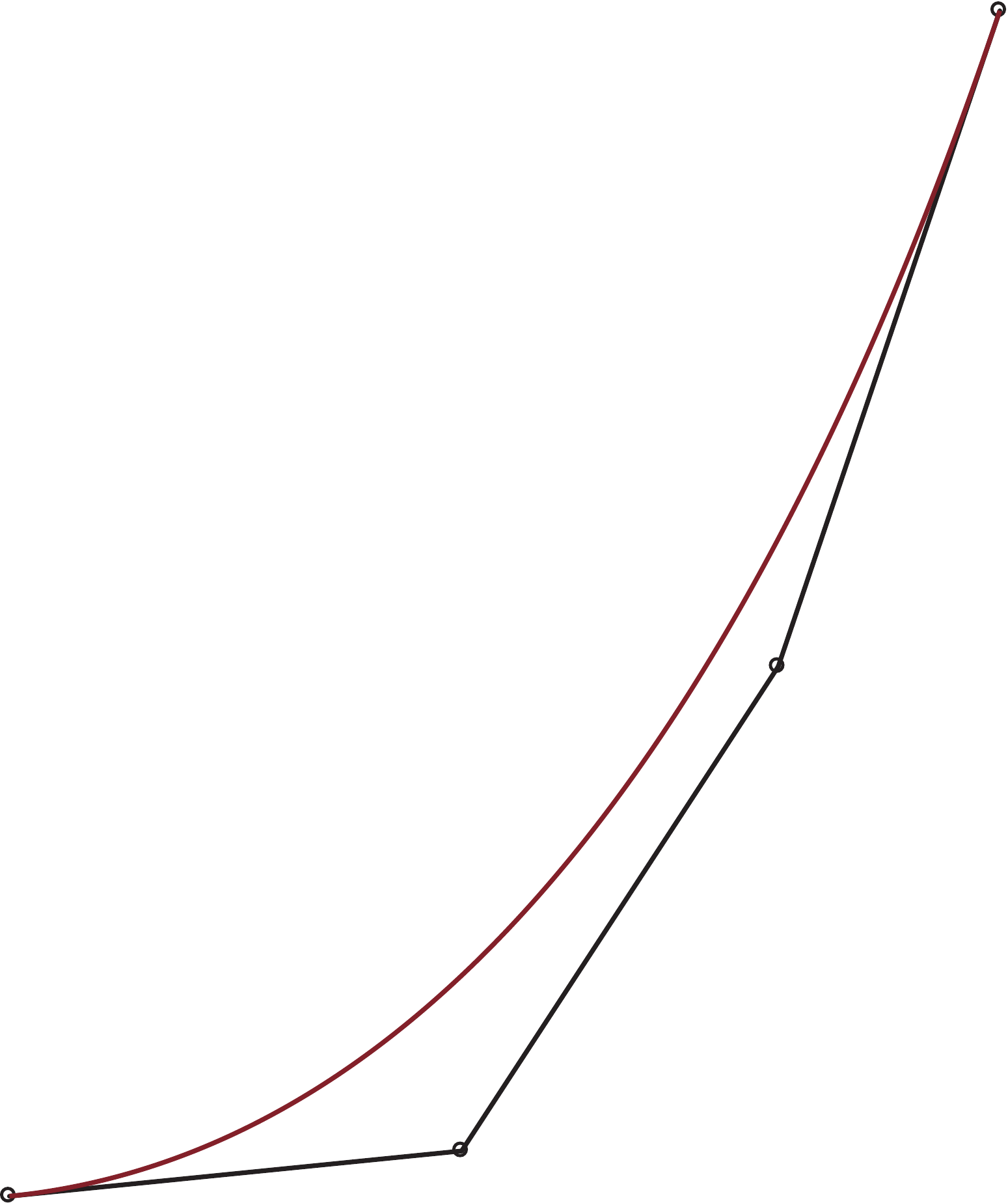}
\end{subfigure}
\hspace{1cm}
\begin{subfigure}{0.3\textwidth}
\includegraphics[width=\textwidth]{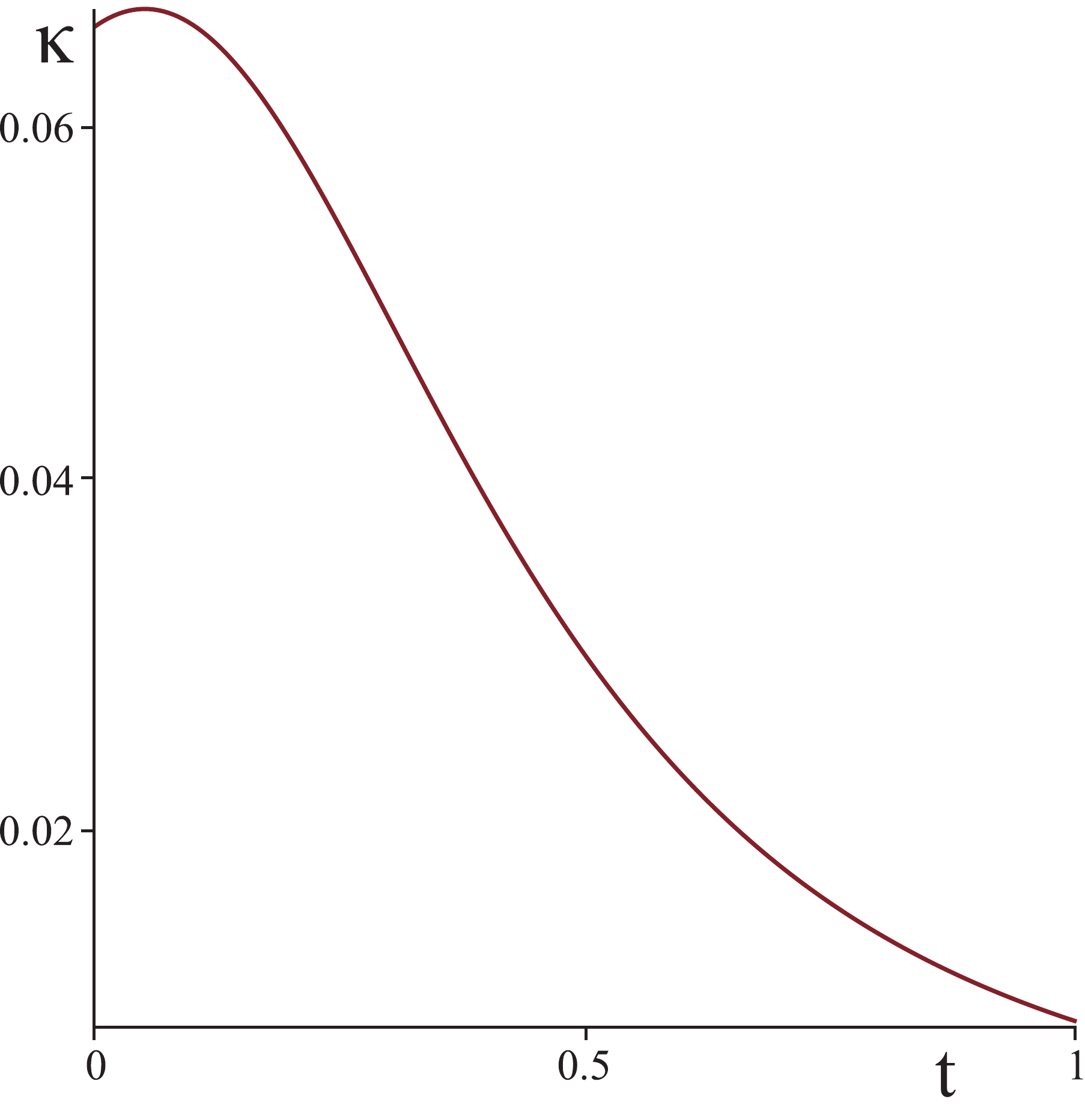}
\end{subfigure}
\caption{Cubic B\'ezier curve generated by $M_6$ and $\mathbf{w}_6$ 
and the graph of its curvature}\label{eje4}
\end{figure}
\end{example}

\begin{example}
Let $M_{7}=\left(\begin{smallmatrix}1 & 0 \\ 1 & 1\end{smallmatrix}\right)$
and $\mathbf{w}_7=\left(\begin{smallmatrix}1\\-1\end{smallmatrix}\right)$. 
Then the cubic B\'ezier curve generated by $M_7$ and $\mathbf{w}_7$ does not have monotonic curvature. See Figure \ref{eje5}.

\begin{figure}[h]
\centering
\begin{subfigure}{0.4\textwidth}
\includegraphics[width=\textwidth]{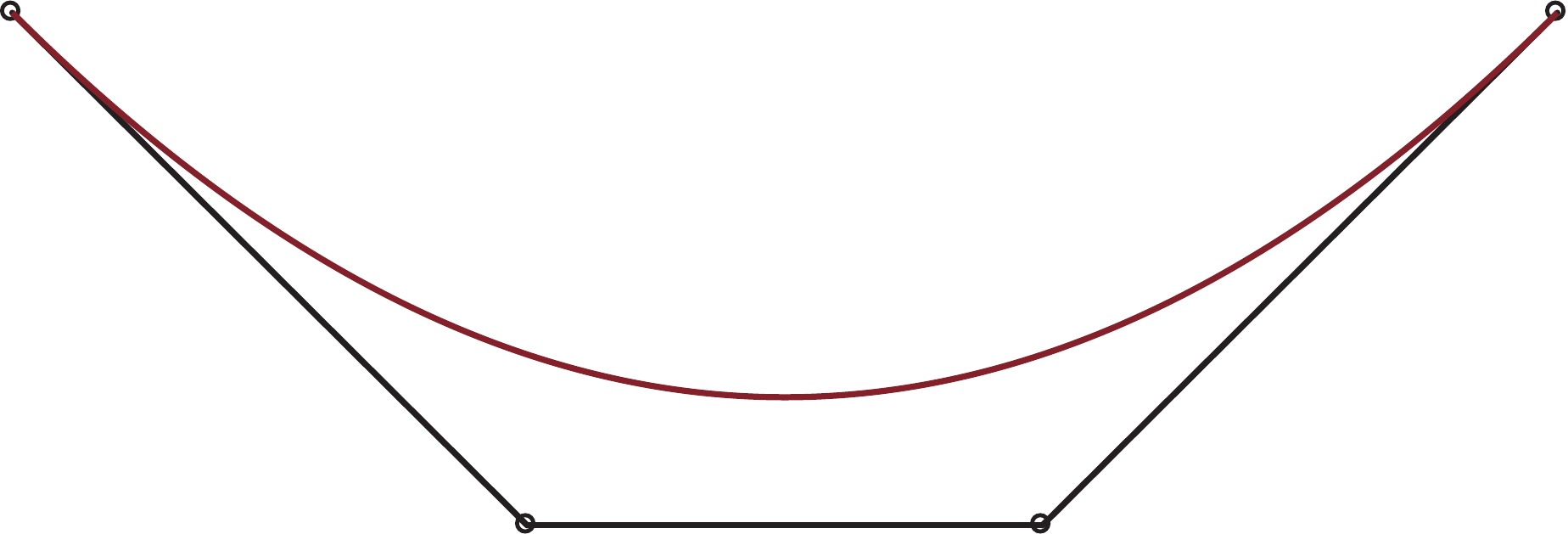}
\end{subfigure}
\hspace{1cm}
\begin{subfigure}{0.3\textwidth}
\includegraphics[width=\textwidth]{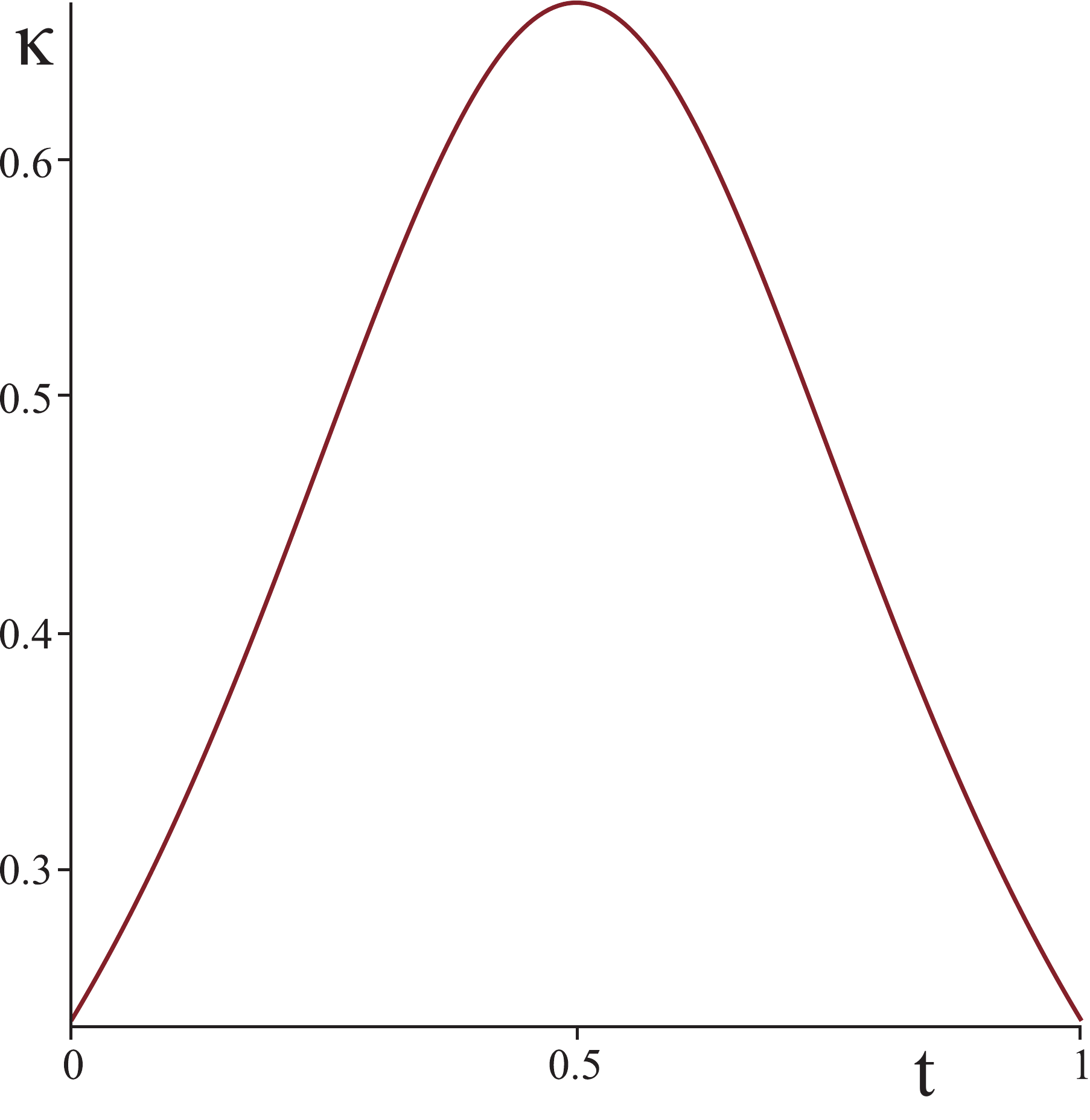}
\end{subfigure}
\caption{Cubic B\'ezier curve generated by $M_7$ and $\mathbf{w}_7$ and the graph of its curvature}\label{eje5}
\end{figure}
\end{example}

We show an example under the conditions of Theorem \ref{tma:Jordan}.

\begin{example}
Let $M_{8}=\left(\begin{smallmatrix}1 & 0 \\ 1 & 1\end{smallmatrix}\right)$
and $\mathbf{w}_8=\left(\begin{smallmatrix}3\\1\end{smallmatrix}\right)$. 
Then the cubic B\'ezier curve generated by $M_8$ and $\mathbf{w}_8$ has
monotonic curvature. See Figure \ref{eje6}.

\begin{figure}[h]
\centering
\begin{subfigure}{0.25\textwidth}
\includegraphics[width=\textwidth]{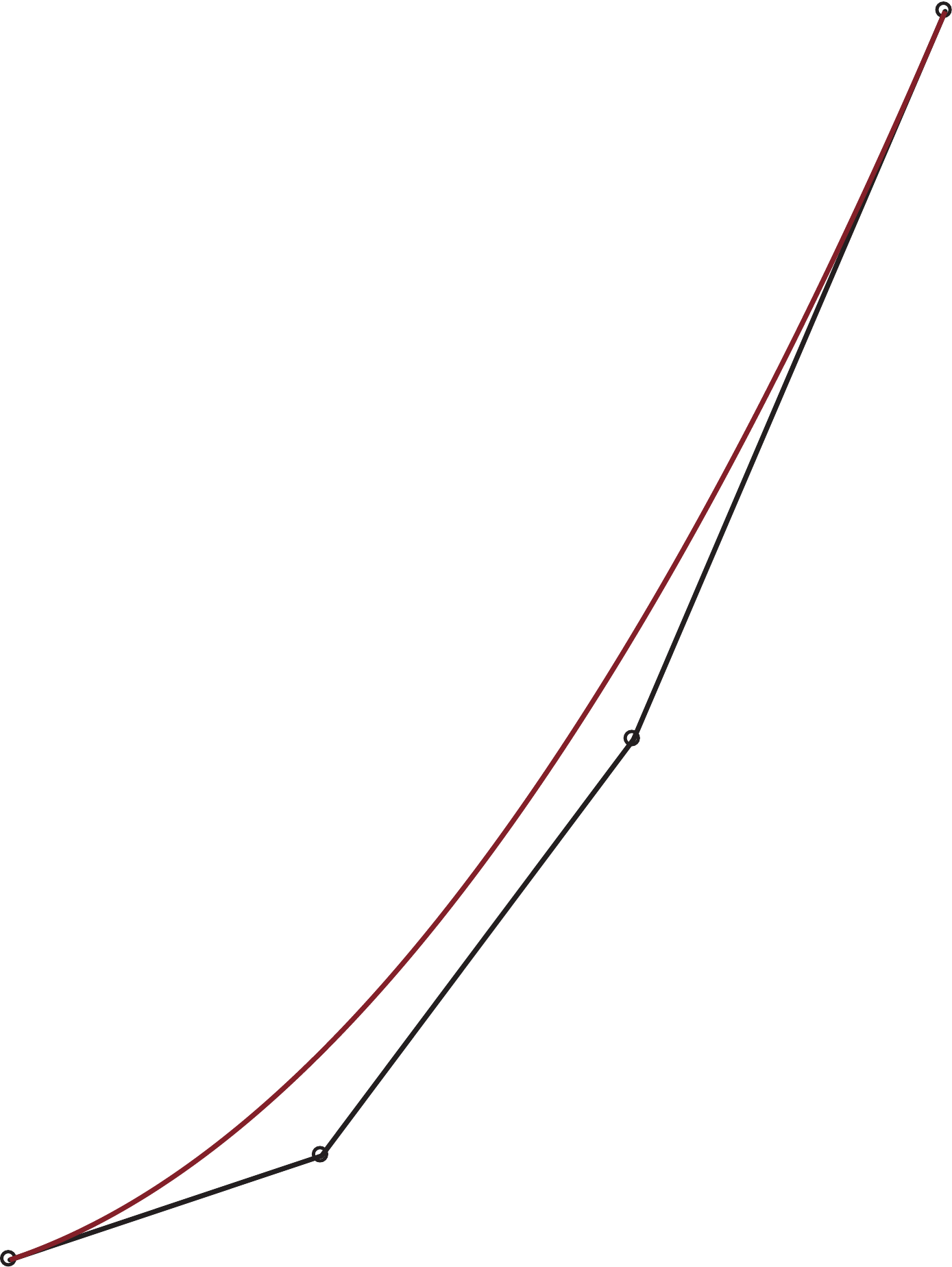}
\end{subfigure}
\hspace{1cm}
\begin{subfigure}{0.3\textwidth}
\includegraphics[width=\textwidth]{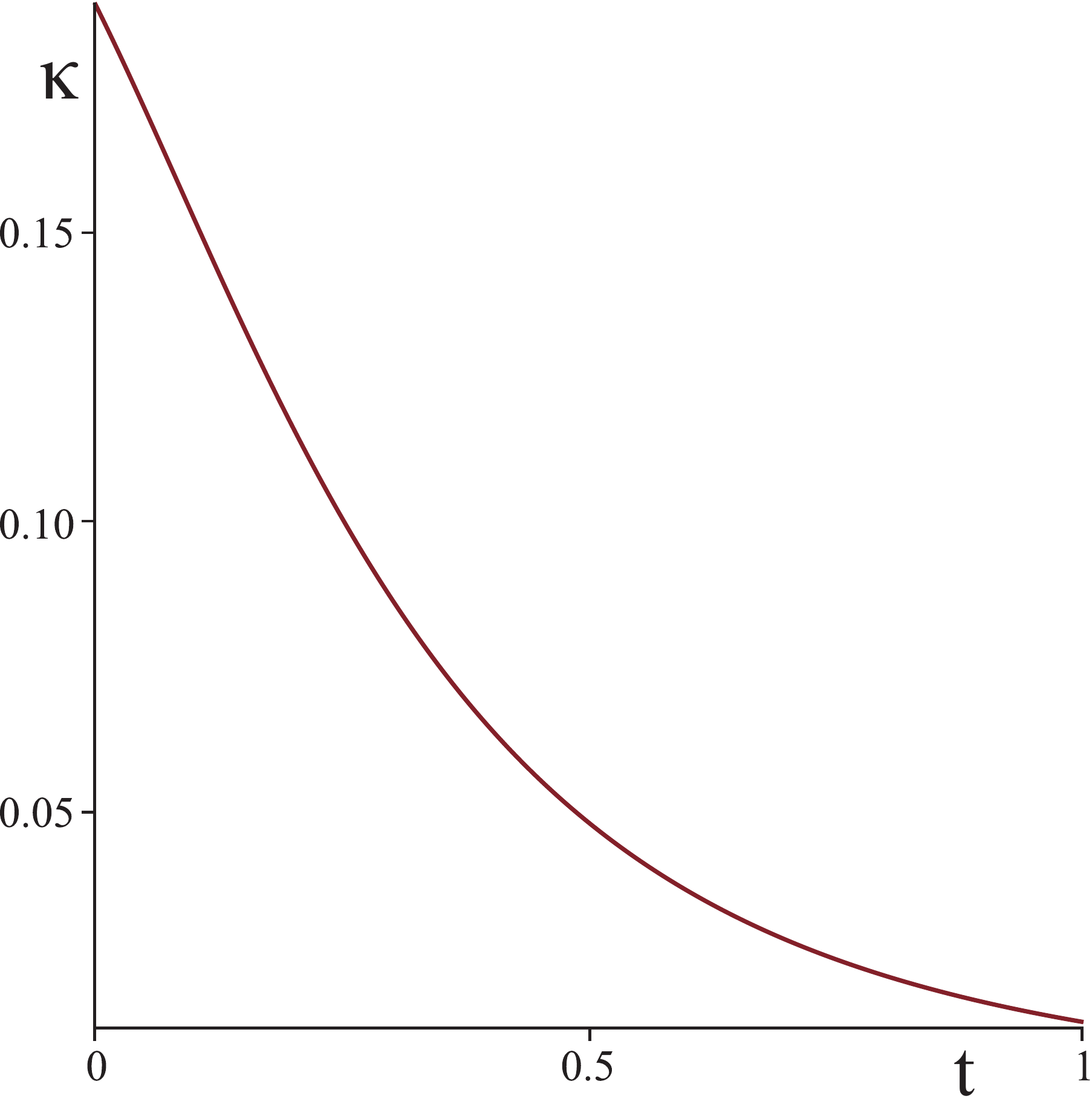}
\end{subfigure}
\caption{Cubic B\'ezier curve generated by $M_8$ and $\mathbf{w}_8$
and the graph of its curvature}\label{eje6}
\end{figure}
\end{example}

% No cumplen

Finally, we show two examples of curves with monotonic curvature for 
which the conditions of this theorem do not hold:

\begin{example}
Consider the matrix $M_{9}=\left(\begin{smallmatrix} 1/2 & -2
\\0 & 1/2\end{smallmatrix}\right)$ with eigenvalue
$\sigma=1/2<1$, with the Jordan basis
$\mathbf{v}_1=\left(\begin{smallmatrix} 1 \\
0\end{smallmatrix}\right)$ and
$\mathbf{v}_2=\left(\begin{smallmatrix} 0 \\
-1/2\end{smallmatrix}\right)$. Then the cubic B\'ezier curve
generated by $M_{9}$ and
$\mathbf{w}_{9}=\tfrac{3}{2}\mathbf{v}_1+2\mathbf{v}_2$ has monotonic
curvature. See Figure \ref{eje14}.

\begin{figure}[h]
\centering
\begin{subfigure}{0.4\textwidth}
\includegraphics[width=\textwidth]{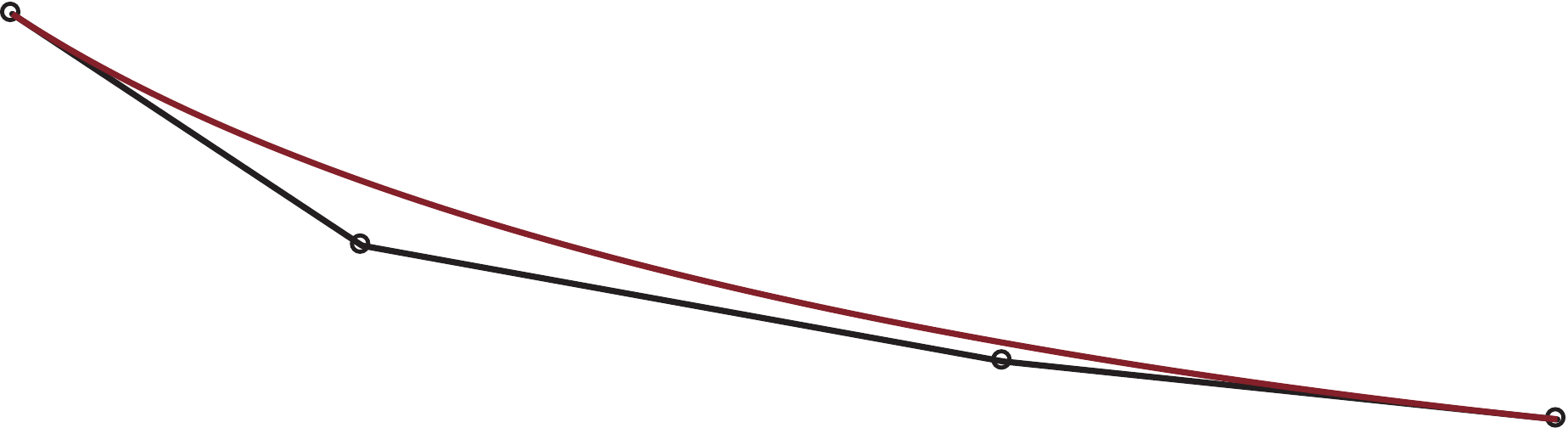}
\end{subfigure}
\hspace{1cm}
\begin{subfigure}{0.3\textwidth}
\includegraphics[width=\textwidth]{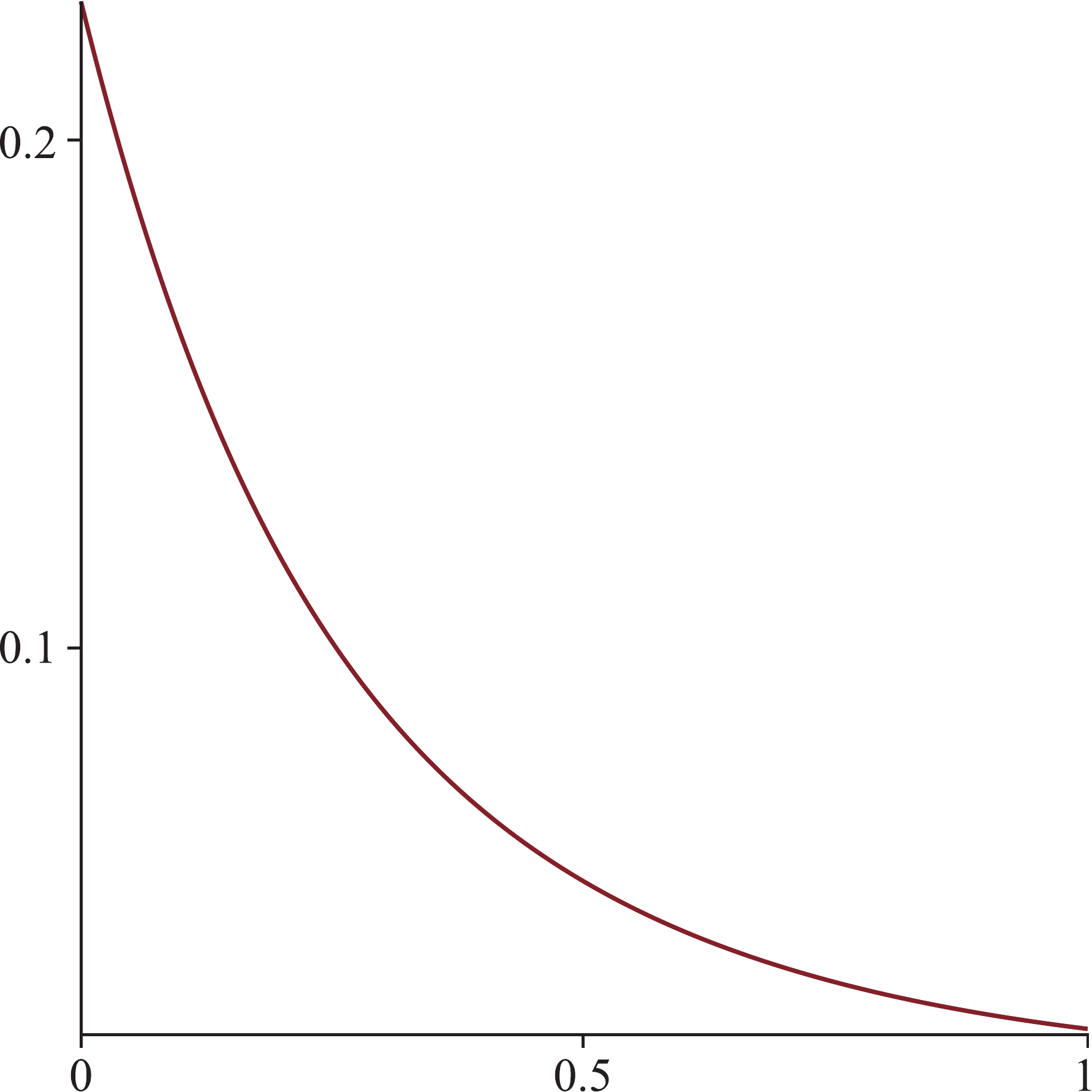}
\end{subfigure}
\caption{Cubic B\'ezier curve generated by $M_9$ and 
$\mathbf{w}_9$
and the graph of its curvature}\label{eje14}
\end{figure}
\end{example}

\begin{example}
Consider the matrix $M_{10}=\left(\begin{smallmatrix} 3/2 & -2
\\0 & 3/2\end{smallmatrix}\right)$ with eigenvalue
$\sigma=3/2$, with the Jordan basis
$\mathbf{v}_1=\left(\begin{smallmatrix} 1 \\
0\end{smallmatrix}\right)$ and
$\mathbf{v}_2=\left(\begin{smallmatrix} 0 \\
-1/2\end{smallmatrix}\right)$. Then the B\'ezier curve of degree $4$
generated by $M_{10}$ and
$\mathbf{w}_{10}=\tfrac{7}{2}\mathbf{v}_1-\tfrac{3}{2}\mathbf{v}_2$ has
monotonic curvature. See Figure \ref{eje15}.

\begin{figure}[h]
\centering
\begin{subfigure}{0.4\textwidth}
\includegraphics[width=\textwidth]{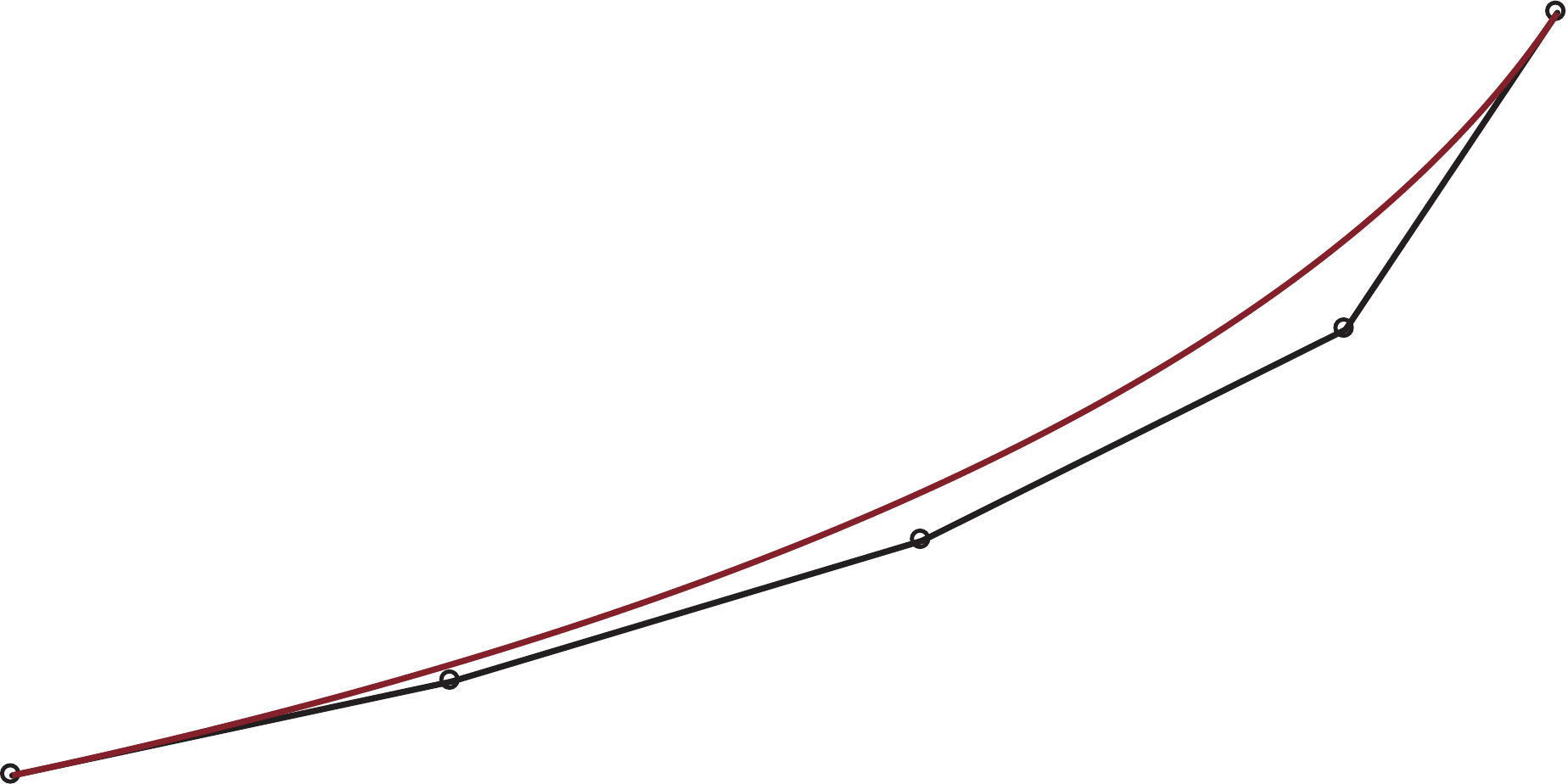}
\end{subfigure}
\hspace{1cm}
\begin{subfigure}{0.3\textwidth}
\includegraphics[width=\textwidth]{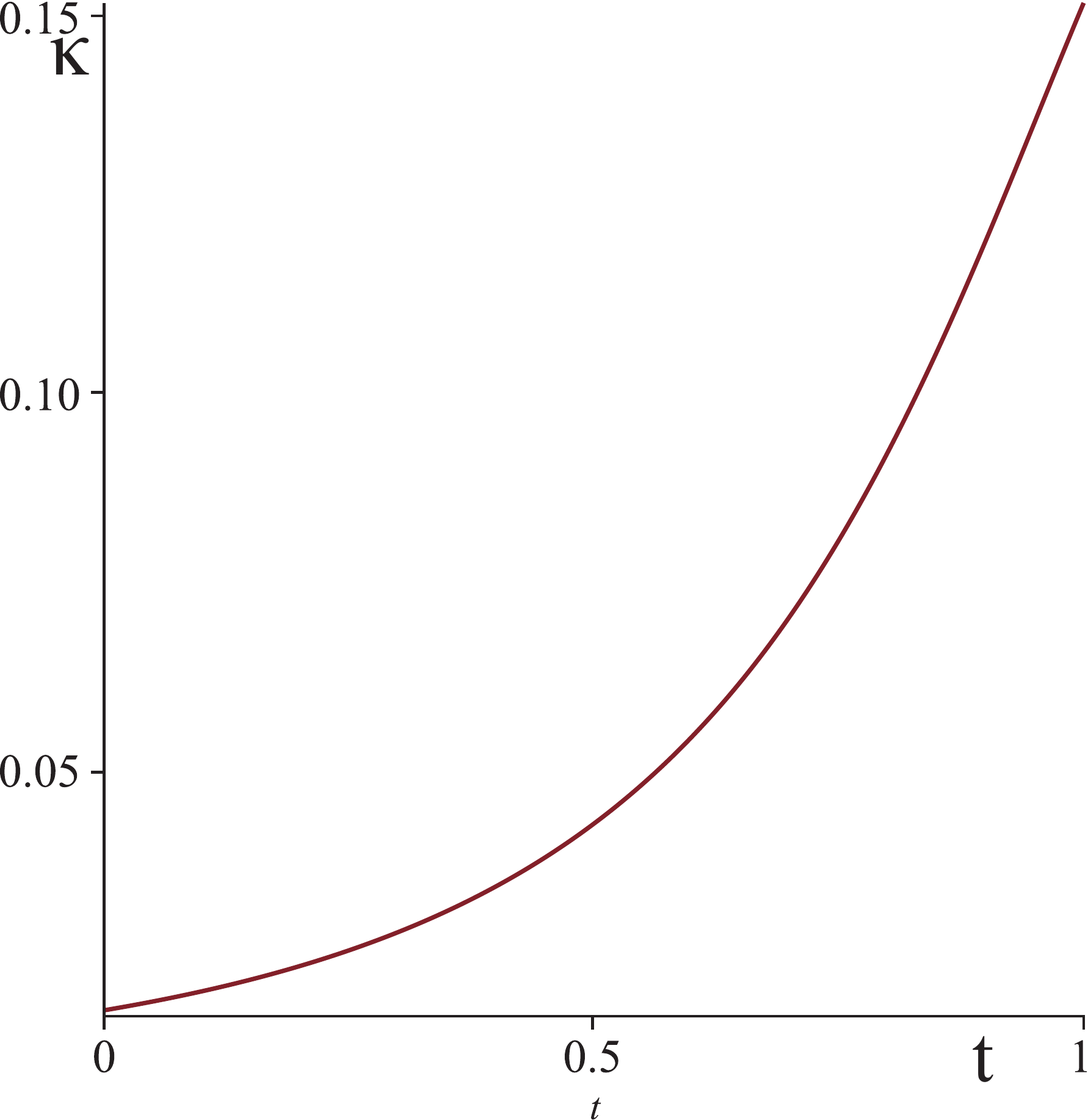}
\end{subfigure}
\caption{Quartic B\'ezier curve generated by $M_{10}$ and 
$\mathbf{w}_{10}$
and the graph of its curvature}\label{eje15}
\end{figure}
\end{example}

Consider now a diagonalizable matrix $M$ with complex eigenvalues $\sigma_1=\sigma$ and $\sigma_2=\overline{\sigma}$. In this case the curvature (\ref{formula_curvatura}) can be written as
\begin{equation}\label{curvatura_compleja}
\kappa(t)=\kappa(0)\dfrac{\vert\sigma(t)\vert^{2(n-2)}\Vert\mathbf{w}\Vert^3}{\Vert T^{n-1}\mathbf{w}\Vert^3},
\end{equation}
where we are writing $\sigma(t)=\sigma_1(t)$ and $\overline{\sigma}(t)=\sigma_2(t)$.

\

As it was mentioned above, in \cite{mineur}, Mineur et al.  define
{\em typical curves} as those planar B\'ezier curves of degree $n$ for
which the edges of its control polygon are obtained iteratively by a
rotation of angle $\varphi$, $|\varphi|<\pi/2$, and a scaling of
factor $h>0$ of the previous edge.  They also give an explicit
expression of the curvature that can be easily recovered using
(\ref{curvatura_compleja}).

\begin{theorem}[Mineur-Lichah-Castelain-Giaume]
Let $h>0$ and $|\varphi|<\pi/2$, $\varphi\neq 0$ and let
\[
M=h\begin{pmatrix}\cos\varphi & -\sin\varphi\\ \sin\varphi & \cos\varphi\end{pmatrix}.
\]
If $h>1/\cos\varphi$ or $0<h<\cos\varphi$ then the B\'ezier curve of degree $n\ge 2$ generated by $M$ and any vector $\mathbf{w}\neq \mathbf{0}$ has monotonic curvature.
\end{theorem}

\

\noindent {\it Proof}. The matrix $M$ has eigenvector $\mathbf{v}=\left(\begin{smallmatrix}
1\\i\end{smallmatrix}\right)$ with eigenvalue $\sigma=he^{-i\varphi}$ so $\overline{\sigma}-\sigma=(2h\sin\varphi)i$ and $\det(\mathbf{v},\overline{\mathbf{v}})=-2i$. If $\mathbf{w}=\mu\,\mathbf{v}+\overline{\mu}\,\overline{\mathbf{v}}$ then, $\Vert \mathbf{w}\Vert=2\vert\mu\vert$, and thus by Theorem \ref{tma curvatura}
\[
\kappa(0)=\dfrac{n-1}{n}\dfrac{4|\mu|^2h\sin\varphi}{\|\mathbf{w}\|^{3}}=\dfrac{n-1}{n}\dfrac{h\sin\varphi}{\|\mathbf{w}\|}.
\]
Now, due to the fact that the hermitian product $\mathbf{v}\cdot\overline{\mathbf{v}}=0$ (analogously to the example of a real diagonalizable symmetric matrix $M$), $\Vert T^{n-1}\mathbf{w}\Vert=|\sigma(t)|^{n-1}\Vert\mathbf{w}\Vert$ and substituting in (\ref{curvatura_compleja})
\[
\kappa(t)=\dfrac{\kappa(0)}{\vert\sigma(t)\vert^{n+1}}=\dfrac{n-1}{n}\dfrac{h\sin\varphi}{\Vert\mathbf{w}\Vert\vert\sigma(t)\vert^{n+1}}.
\]
The derivate takes a simple form,
\[
\kappa'(t)=-\dfrac{n+1}{2}\dfrac{\kappa(0)}{\vert\sigma(t)\vert^{n+3}}\bigl((\sigma+\overline{\sigma}-2)(1-t)+(2\sigma\overline{\sigma}-\sigma-\overline{\sigma})t\bigr).%\vert\sigma(t)\vert^2\bigr)'
\]

From the expression above it is immediately deduced that the curvature is monotonic if and only if the last factor in the parenthesis does not change its sign for $t\in[0,1]$. Since $\sigma=he^{-i\varphi}$
\[
(\sigma+\overline{\sigma}-2)(1-t)+(2\sigma\overline{\sigma}-\sigma-\overline{\sigma})t=2(h\cos\varphi-1)(1-t)+2h(h-\cos\varphi)t,
\]
thus, in this situation, the curvature is monotonic if and only if $h>1/\cos\varphi$ or $0<h<\cos\varphi$, as shown in \cite{mineur}. \hfill $\square$

\begin{example}
Let $h=1.8$, $\varphi=0.925$, so that 
$M_{11}=\left(\begin{smallmatrix}1.083 & -1.438 \\ 1.438 &
1.083\end{smallmatrix}\right)$ and
$\mathbf{w}_{11}=\left(\begin{smallmatrix}0.4\\0.1\end{smallmatrix}\right)$.
Then the B\'ezier curve generated by $M_{11}$ and $\mathbf{w}_{11}$
has monotonic curvature.  For Figure \ref{eje7} we have taken $n=7$.

\begin{figure}[h]
\centering
\begin{subfigure}{0.3\textwidth}
\includegraphics[width=\textwidth]{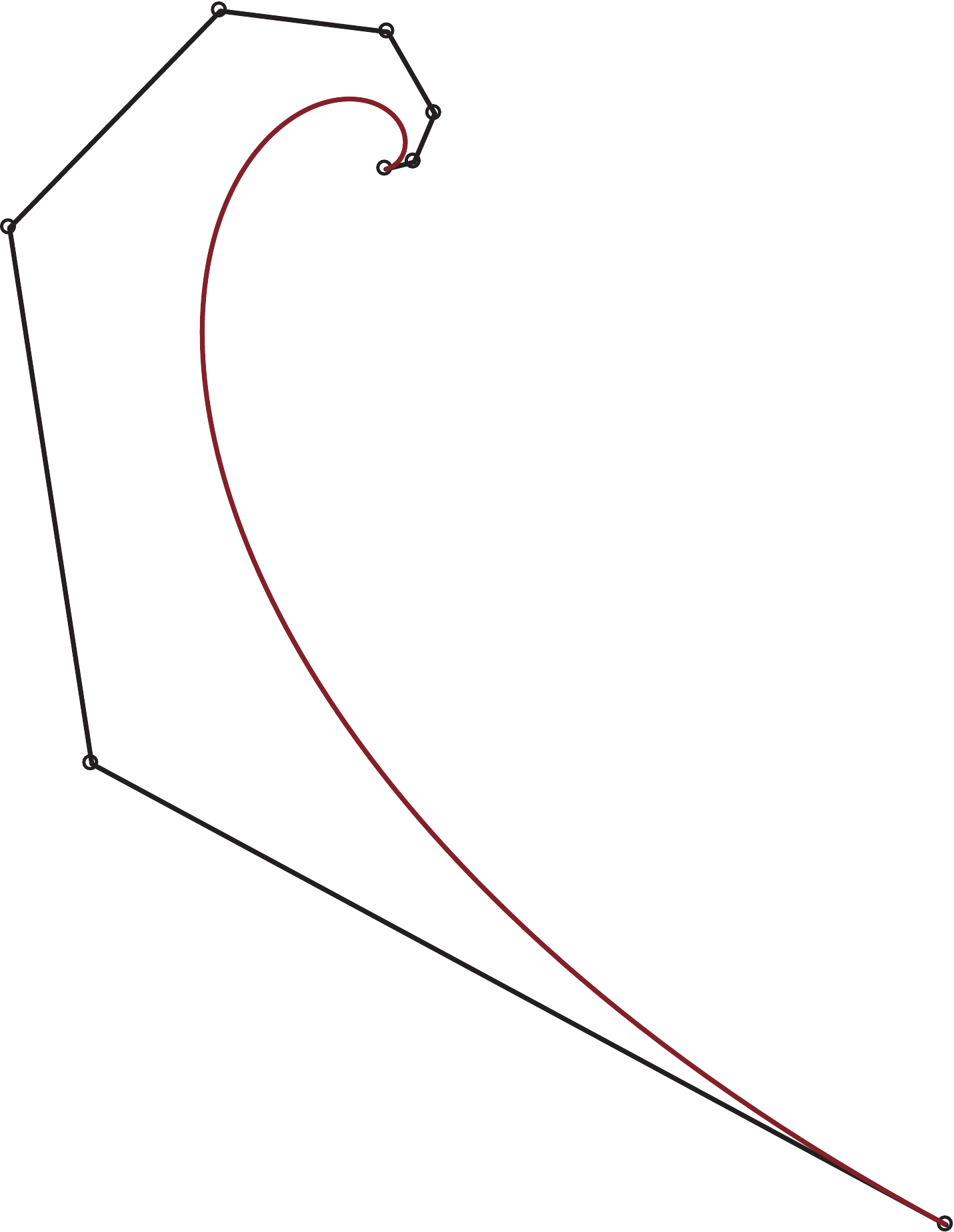}
\end{subfigure}
\hspace{1cm}
\begin{subfigure}{0.3\textwidth}
\includegraphics[width=\textwidth]{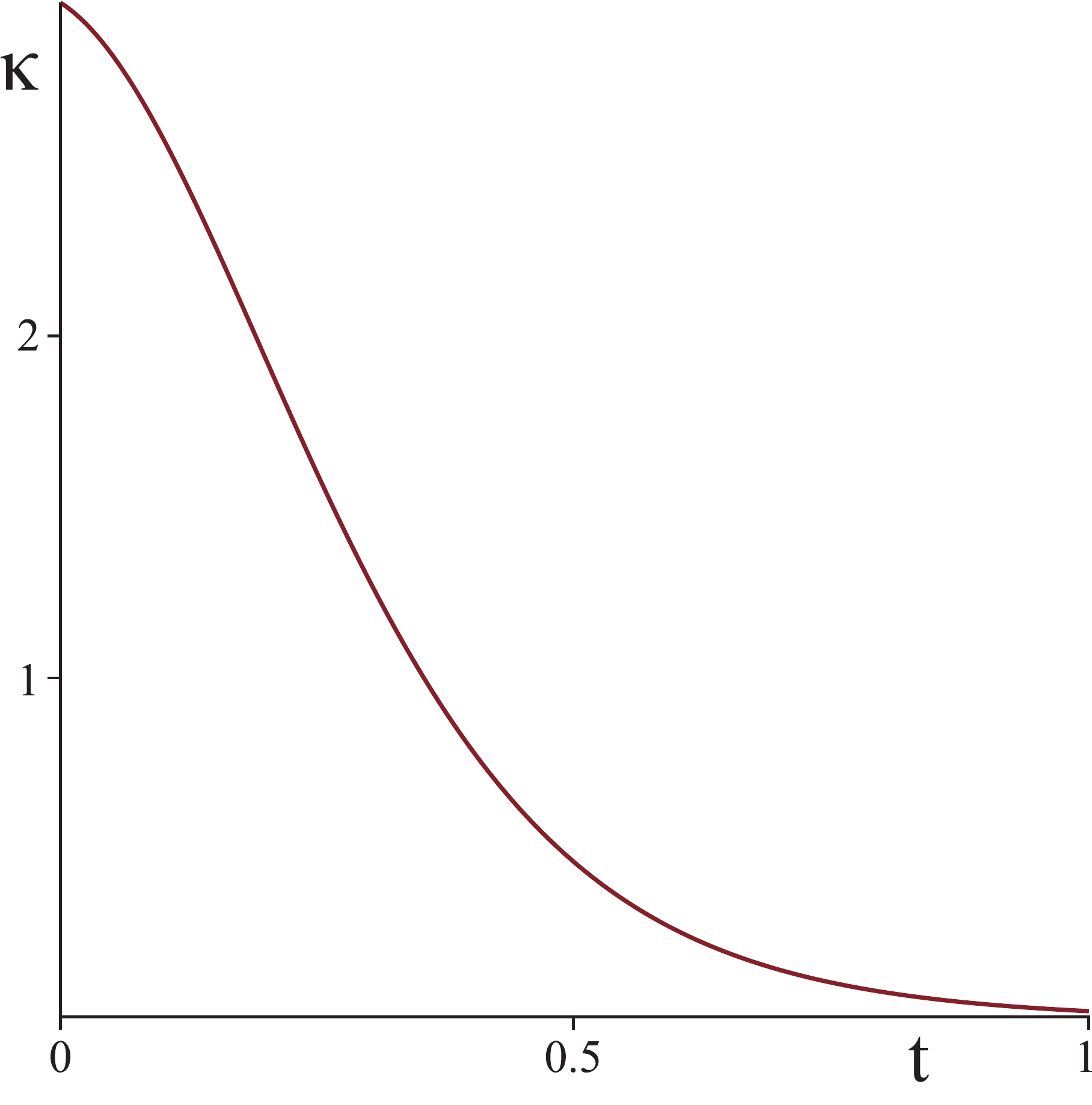}
\end{subfigure}
\caption{B\'ezier curve of degree 7 generated by $M_{11}$ and 
$\mathbf{w}_{11}$ and the graph of its curvature}\label{eje7}
\end{figure}
\end{example}

\begin{example}
Let $h=1.2$, $\varphi=0.925$, so that 
$M_{12}=\left(\begin{smallmatrix}0.722 & -0.958 \\ 0.958 &
0.722\end{smallmatrix}\right)$ and
$\mathbf{w}_{12}=\left(\begin{smallmatrix}0.4\\0.1\end{smallmatrix}\right)$.
The B\'ezier curve generated by $M_{12}$ and $\mathbf{w}_{12}$
does not have monotonic curvature.  For Figure \ref{eje8} we have taken $n=7$.

\begin{figure}[h]
\centering
\begin{subfigure}{0.3\textwidth}
\includegraphics[width=\textwidth]{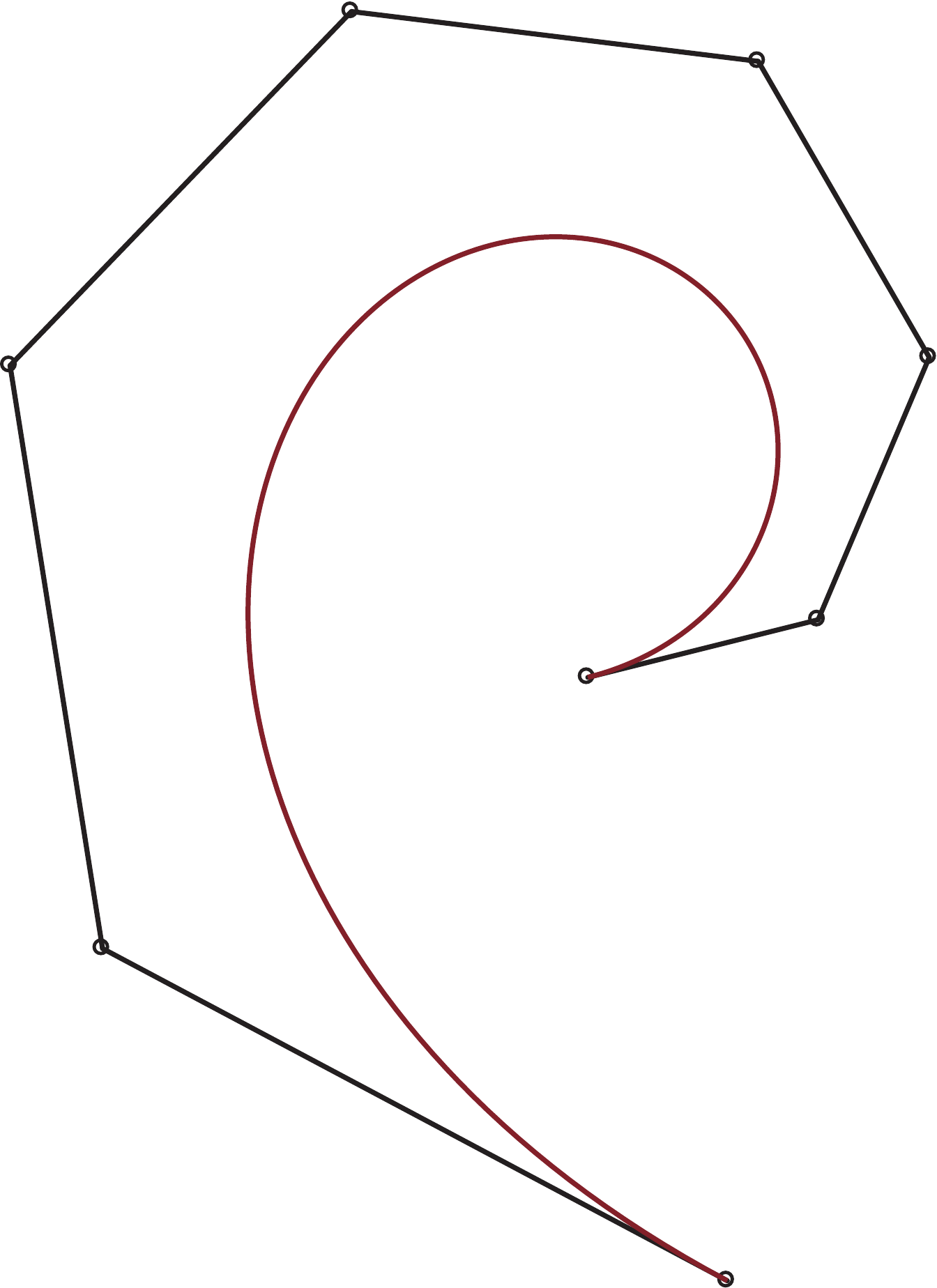}
\end{subfigure}
\hspace{1cm}
\begin{subfigure}{0.3\textwidth}
\includegraphics[width=\textwidth]{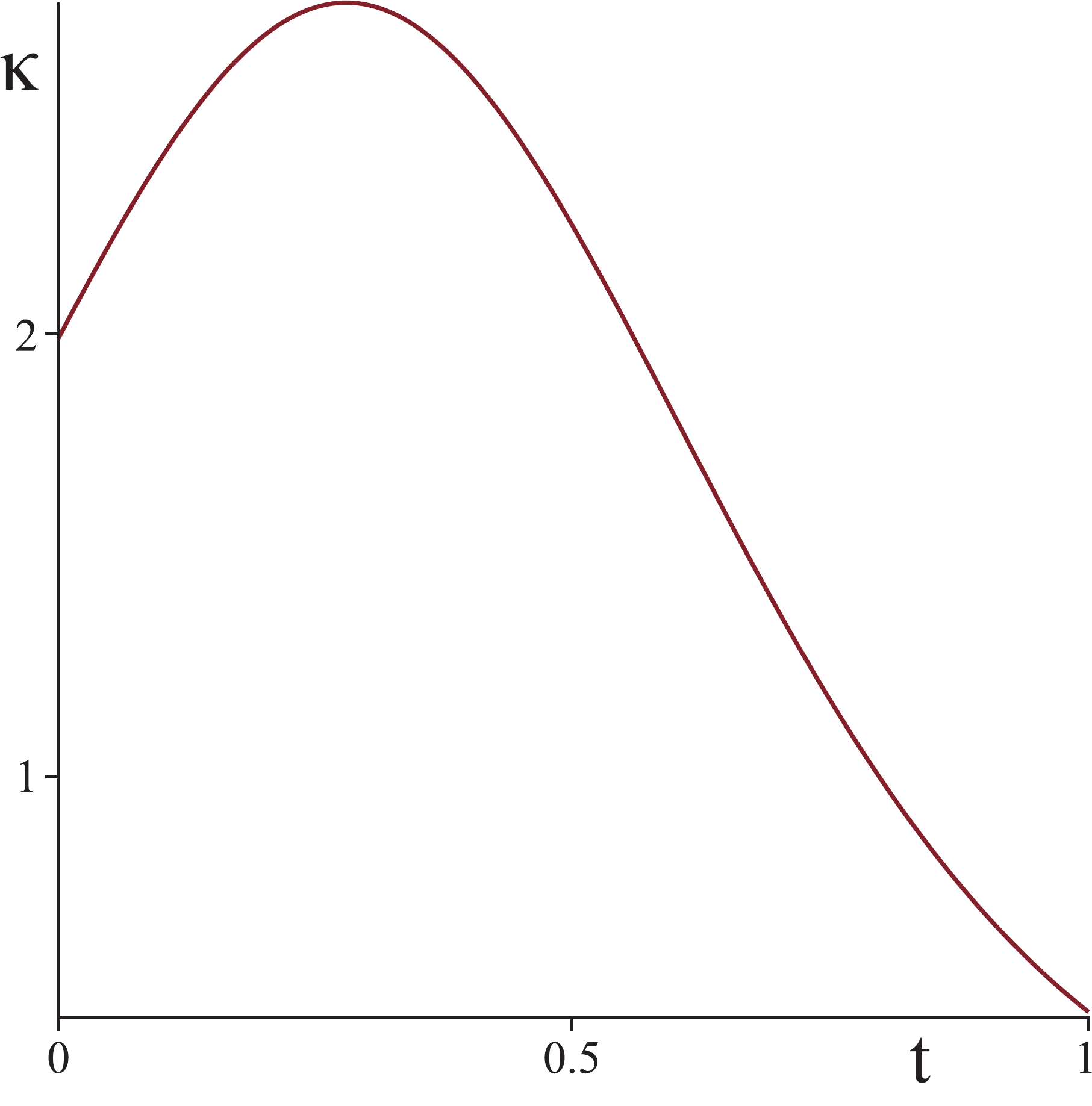}
\end{subfigure}
\caption{B\'ezier curve of degree 7 generated by $M_{12}$ and
$\mathbf{w}_{12}$ and the graph of its curvature}\label{eje8}
\end{figure}
\end{example}

With the explicit expression for the curvature, it is possible to extend the previous result to more general matrices $M$ with complex eigenvalues:

\begin{theorem}\label{t:complejo}
Let $M$ be a diagonalizable matrix with non-real complex eigenvalues $\sigma=he^{i\varphi}$ and $\overline{\sigma}=he^{-i\varphi}$, and corresponding eigenvectors $\mathbf{v}$ and $\overline{\mathbf{v}}$. Let $\gamma$ be the angle formed by the vectors $\operatorname{Re}\mathbf{v}$ and $\operatorname{Im}\mathbf{v}$. If
\begin{equation}\label{exMineur}
\left\vert\cos\gamma\right\vert< \dfrac{h\cos\varphi-1}{\sqrt{(h\cos\varphi-1)^2+9h^2\sin^2\varphi}}
\end{equation}
then the B\'ezier curve of degree $n\ge 2$ generated by $M$ and any vector $\mathbf{w}\neq \mathbf{0}$ has monotonic curvature (decreasing if $\kappa(0)>0$ and increasing if $\kappa(0)<0$).

Also, if
\[
\left\vert\cos\gamma\right\vert< \dfrac{\cos\varphi-h}{\sqrt{(\cos\varphi-h)^2+9\sin^2\varphi}}
\]
then the B\'ezier curve of degree $n\ge 2$ generated by $M$ and any vector $\mathbf{w}\neq \mathbf{0}$ has monotonic curvature (increasing if $\kappa(0)>0$ and decreasing if $\kappa(0)<0$).
\end{theorem}

Observe that when $\operatorname{Re}\mathbf{v}$ and $\operatorname{Im}\mathbf{v}$ are orthogonal (that is $\gamma=\frac{\pi}{2}$ or $\gamma=-\frac{\pi}{2}$) we recover the result of Mineur et al. stated in the previous theorem.

\

\noindent {\it Proof}. The eigenvector $\mathbf{v}$ is chosen so that $\Vert \operatorname{Re}\mathbf{v} \Vert=\Vert \operatorname{Im}\mathbf{v} \Vert$. Write $\mathbf{w}\neq \mathbf{0}$ as $\mathbf{w}=\mu\mathbf{v}+\overline{\mu\mathbf{v}}$ with $\mu\in\mathbb{C}\setminus\{0\}$.

We use the expression of the derivative of the curvature shown in Lemma \ref{l:simplificaderivada} below:
\[
\begin{split}
\kappa'(t)&=\dfrac{2\kappa(0)\Vert\mathbf{v}\Vert^2\Vert \mathbf{w}\Vert ^{3}\vert\mu\vert^2\vert\sigma(t)\vert^{2(2n-4)}}{\Vert T^{n-1}\mathbf{w}
\Vert^5}\\
&\cdot
\Bigl[(n+1)M(t)\bigl(-1+\cos\gamma\operatorname{Im}\left(e^{i\theta}e^{i2(n-1)\varphi(t)}\right)\bigr)\Bigr.\\
&\hspace{0.5cm}\Bigl.+3(n-1)h\cos\gamma\sin\varphi \operatorname{Re}\left(e^{i\theta}e^{i2(n-1)\varphi(t)}\right)\Bigr],
\end{split}
\]
where $M(t)=(h\cos\varphi-1)(1-t)+t(h^2-h\cos\varphi)$, $\sigma(t)=1-t+t\sigma=\vert\sigma(t)\vert e^{i\varphi(t)}$, and $\theta$ is an argument of $\mu^2$. Hence the sign of $\kappa'(t)/\kappa(0)$ is given by the sign of
\[
\begin{split}
S(t)=&(n+1)M(t)\bigl(-1+\cos\gamma \operatorname{Im}\left(e^{i\theta}e^{i2(n-1)\varphi(t)}\right)\bigr)\\
&+3(n-1)h\cos\gamma\sin\varphi\operatorname{Re}\left(e^{i\theta}e^{i2(n-1)\varphi(t)}\right).
\end{split}
\]
Assume that (\ref{exMineur}) holds, then $h\cos\varphi-1>0$ and hence $M(t)>0$ for every $t\in[0,1]$. So $S$ can be bounded by
\[
\begin{split}
S(t)\le &(n+1)M(t)\bigl(-1+\bigl\vert\cos\gamma \operatorname{Im}\left(e^{i\theta}e^{i2(n-1)\varphi(t)}\right)\bigr\vert\bigr)\\
&+3(n-1)h\bigl\vert\cos\gamma\sin\varphi\operatorname{Re}\left(e^{i\theta}e^{i2(n-1)\varphi(t)}\right)\bigr\vert\\
\end{split}
\]
Consider now the terms with common factor $\left\vert\cos\gamma\right\vert$, that is,
\[
(n+1)M(t)\bigl\vert\operatorname{Im}\left(e^{i\theta}e^{i2(n-1)\varphi(t)}\right)\bigr\vert+3(n-1)h\bigl\vert\sin\varphi\operatorname{Re}\left(e^{i\theta}e^{i2(n-1)\varphi(t)}\right)\bigr\vert.
\]
This expression gets its biggest value when $\bigl\vert\operatorname{Im}\left(e^{i\theta}e^{i2(n-1)\varphi(t)}\right)\bigr\vert$ is replaced by
\[
\dfrac{(n+1)M(t)}{\sqrt{(n+1)^2M^2(t)+9(n-1)^2h^2\sin^2\varphi}}
\]
and $\bigl\vert\operatorname{Re}\left(e^{i\theta}e^{i2(n-1)\varphi(t)}\right)\bigr\vert$ is replaced by
\[
\dfrac{3(n-1)h\left\vert\sin\varphi\right\vert}{\sqrt{(n+1)^2M^2(t)+9(n-1)^2h^2\sin^2\varphi}}.
\]
Hence $S$ can be bounded by
\[
S(t)\le -(n+1)M(t)+\left\vert\cos\gamma\right\vert\sqrt{(n+1)^2M^2(t)+9(n-1)^2h^2\sin^2\varphi},
\]
and by (\ref{exMineur}) we obtain,
\[
S(t)\le -(n+1)M(t)+\dfrac{(h\cos\varphi-1)\sqrt{(n+1)^2M^2(t)+9(n-1)^2h^2\sin^2\varphi}}{\sqrt{(h\cos\varphi-1)^2+9h^2\sin^2\varphi}}
\]
Now, since for any $n\ge 2$ we have $0< (n-1)/(n+1)<1$, 
\[
\dfrac{h\cos\varphi-1}{\sqrt{(h\cos\varphi-1)^2+9h^2\sin^2\varphi}}< \dfrac{(n+1)(h\cos\varphi-1)}{\sqrt{(n+1)^2(h\cos\varphi-1)^2+9(n-1)^2h^2\sin^2\varphi}}
\]
and since $M(t)\ge h\cos\varphi-1>0$ for any $t\in[0,1]$, this last quotient is less or equal than
\[
\dfrac{(n+1)M(t)}{\sqrt{(n+1)^2M^2(t)+9(n-1)^2h^2\sin^2\varphi}}.
\]
Therefore, we obtain that
\[
\dfrac{(h\cos\varphi-1)}{\sqrt{(h\cos\varphi-1)^2+9h^2\sin^2\varphi}}<\dfrac{(n+1)M(t)}{\sqrt{(n+1)^2M^2(t)+9(n-1)^2h^2\sin^2\varphi}}
\]
and thus
\[
S(t)< -(n+1)M(t)+(n+1)M(t)=0,
\]
so $\kappa'(t)/\kappa(0)<0$ as claimed.

For the other case the reasoning is similar but using that $M(t)<0$ for every $t\in[0,1]$ and, that the minimum of $-M(t)$ is attained at $t=1$ and given by
$-M(1)=h(\cos\varphi-h)$. \hfill $\square$

\

The proof of Theorem \ref{t:complejo} gives a condition to get
B\'ezier curves of a given degree with monotonic curvature.
Concretely:

\begin{corollary}
Let $n\ge 2$ be a given integer and $M$, $\sigma=he^{i\varphi}$, and $\mathbf{v}$ as in Theorem \ref{t:complejo}. If 
\begin{equation}\label{exMineur_grado}
\left\vert\cos\gamma\right\vert\le \dfrac{(n+1)(h\cos\varphi-1)}{\sqrt{(n+1)^2(h\cos\varphi-1)^2+9(n-1)^2h^2\sin^2\varphi}} 
\end{equation}
or
\[
\left\vert\cos\gamma\right\vert\le \dfrac{(n+1)(\cos\varphi-h)}{\sqrt{(n+1)^2(\cos\varphi-h)^2+9(n-1)^2\sin^2\varphi}}
\]
then the B\'ezier curve of degree $n$ generated by $M$ and $\mathbf{w}\neq\mathbf{0}$ has monotonic curvature.
\end{corollary}

\

The next result is a technical lemma to obtain the expression for the derivative of the curvature used in Theorem \ref{t:complejo}:

\begin{lemma}\label{l:simplificaderivada}
Let $M$ be a diagonalizable matrix with non-real complex eigenvalue $\sigma=he^{i\varphi}$ and, choose its eigenvector $\mathbf{v}$ so that $\left\Vert\operatorname{Re}\mathbf{v}\right\Vert=\left\Vert\operatorname{Im}\mathbf{v}\right\Vert$. Let $\mathbf{w}\neq \mathbf{0}$ be written as $\mathbf{w}=\mu\mathbf{v}+\overline{\mu\mathbf{v}}$.

Then the derivative of the curvature of the B\'ezier curve of degree $n\ge 2$ generated by $M$ and $\mathbf{w}$ is
\[
\begin{split}
\kappa'(t)&=\dfrac{2\kappa(0)\Vert\mathbf{v}\Vert^2\Vert \mathbf{w}\Vert ^{3}\vert\mu\vert^2\vert\sigma(t)\vert^{2(2n-4)}}{\Vert T^{n-1}\mathbf{w}
	\Vert^5}\\
&\cdot
\Bigl[(n+1)\bigl(-1+\cos\gamma\operatorname{Im}\left(z(t)\right)\bigr)\bigl((h\cos\varphi-1)(1-t)+t(h^2-h\cos\varphi)\bigr)\Bigr.\\
&\hspace{0.5cm}\Bigl.+3(n-1)h\cos\gamma\sin\varphi \operatorname{Re}\left(z(t)\right)\Bigr].
\end{split}
\] 
where $\gamma$ is the angle formed by the vectors $\operatorname{Re}\mathbf{v}$ and $\operatorname{Im}\mathbf{v}$, $\mu^2=\vert\mu\vert^2e^{i\theta}$ and $\sigma(t)=\vert\sigma(t)\vert e^{i\varphi(t)}=1-t+\sigma t$ is an eigenvalue of $T=(1-t)\mathbb{I}+tM$ and $z(t)=e^{i\theta}e^{i2(n-1)\varphi(t)}$.
\end{lemma}

\

\noindent {\it Proof}. Consider the derivative given by (\ref{dcurvature}) with $\sigma_1(t)=\sigma(t)=1-t+\sigma t$ and $\sigma_2(t)=\overline{\sigma}(t)=1-t+\overline{\sigma}t$. The norm of
$T^{n-1}\mathbf{w}$ is calculated by means of the hermitian product,
\[
\Vert T^{n-1}\mathbf{w}\Vert^2=2\Vert\mathbf{v}\Vert^2\bigl(\vert\mu\vert^2\vert\sigma(t)\vert^{2(n-1)}-\cos\gamma\operatorname{Im}\left(\mu^2\sigma(t)^{2(n-1)}\right)\bigr),
\]
since $\mathbf{v}\cdot\overline{\mathbf{v}}=i\cos\gamma \Vert \mathbf{v}\Vert^2$ because of the choice of $\mathbf{v}$ with $\left\Vert\operatorname{Re}\mathbf{v}\right\Vert=\left\Vert\operatorname{Im}\mathbf{v}\right\Vert$.

Substituting in (\ref{dcurvature}) we obtain
\[
\begin{split}
\kappa'(t)&=\dfrac{\kappa(0)\Vert \mathbf{w}\Vert ^{3}\vert\sigma(t)\vert^{2(n-3)}}{\Vert T^{n-1}\mathbf{w}
\Vert^5}\\
&\cdot
\Bigl[2(n-2)\Vert\mathbf{v}\Vert^2\left(\vert\mu\vert^2\vert\sigma(t)\vert^{2(n-1)}-\cos\gamma\operatorname{Im}\left(\mu^2\sigma(t)^{2(n-1)}\right)\right)\bigl(\vert\sigma(t)\vert^2\bigr)'\Bigr.\\
&\hspace{0.5cm}\Bigl.-3\vert\sigma(t)\vert^2 \Vert \mathbf{v}\Vert^2\bigl(\vert\mu\vert^2\bigl(\vert\sigma(t)\vert^{2(n-1)}\bigr)'-\cos\gamma \operatorname{Im}\left(\mu^2\left(\sigma(t)^{2(n-1)}\right)'\right)\Bigr].
\end{split}
\]
Taking the derivatives of $\vert\sigma(t)\vert^{2(n-1)}$ and of $\sigma(t)^{2(n-1)}$ as powers of $\vert\sigma(t)\vert^2$ and $\sigma^2(t)$ respectively, and taking common factors $\Vert\mathbf{v}\Vert^2$, $\vert\mu\vert^2$ and $\vert\sigma(t)\vert^{2(n-1)}$
we get
\[
\begin{split}
\kappa'(t)&=\dfrac{\kappa(0)\Vert\mathbf{v}\Vert^2\Vert \mathbf{w}\Vert ^{3}\vert\mu\vert^2\vert\sigma(t)\vert^{2(2n-4)}}{\Vert T^{n-1}\mathbf{w}
\Vert^5}\\
&\cdot
\Bigl[2(n-2)\left(1-\cos\gamma\operatorname{Im}\left(e^{i\theta}e^{i2(n-1)\varphi(t)}\right)\right)\bigl(\vert\sigma(t)\vert^2\bigr)'\Bigr.\\
&\hspace{0.5cm}\Bigl.-3(n-1)\left(\bigl(\vert\sigma(t)\vert^{2}\bigr)'-\cos\gamma \operatorname{Im}\left(e^{i\theta}e^{i2(n-2)\varphi(t)}\left(\sigma^2(t)\right)'\right)\right)\Bigr].
\end{split}
\]
Thus, grouping the terms multiplied by the derivative of $\vert\sigma(t)\vert^2$,
\[
\begin{split}
\kappa'(t)&=\dfrac{\kappa(0)\Vert\mathbf{v}\Vert^2\Vert \mathbf{w}\Vert ^{3}\vert\mu\vert^2\vert\sigma(t)\vert^{2(2n-4)}}{\Vert T^{n-1}\mathbf{w}
\Vert^5}\\
&\cdot
\Bigl[\bigl(-(n+1)-2(n-2)\cos\gamma\operatorname{Im}\left(e^{i\theta}e^{i2(n-1)\varphi(t)}\right)\bigr)\bigl(\vert\sigma(t)\vert^2\bigr)'\Bigr.\\
&\hspace{0.5cm}\Bigl.+3(n-1)\cos\gamma \operatorname{Im}\left(e^{i\theta}e^{i2(n-2)\varphi(t)}\left(\sigma^2(t)\right)'\right)\Bigr].
\end{split}
\]
Now
\[
\begin{split}
\operatorname{Im}\left(e^{i\theta}e^{i2(n-2)\varphi(t)}\left(\sigma^2(t)\right)'\right)=&
\operatorname{Im}\left(e^{i\theta}e^{i2(n-1)\varphi(t)}\right)\operatorname{Re}\left(e^{-i2\varphi(t)}\left(\sigma^2(t)\right)'\right)\\
&+\operatorname{Re}\left(e^{i\theta}e^{i2(n-1)\varphi(t)}\right)\operatorname{Im}\left(e^{-i2\varphi(t)}\left(\sigma^2(t)\right)'\right),
\end{split}
\]
and since $\sigma(t)=1-t+t\sigma$, and $\overline{\sigma}(t)=1-t+t\overline{\sigma}=\vert\sigma(t)\vert e^{-i\varphi(t)}$,
\[
\begin{split}
e^{-i2\varphi(t)}\left(\sigma^2(t)\right)'&=\dfrac{\left(\overline{\sigma}(t)\right)^2}{\vert\sigma(t)\vert^2}2\sigma(t)(\sigma-1)=2\overline{\sigma}(t)(\sigma-1)\\
&=2\bigl((\sigma-1)(1-t)+t(\vert\sigma\vert^2-\overline{\sigma})\bigr).
\end{split}
\]
Writing $\sigma=he^{i\varphi}$,
\[
\begin{array}{l}
\operatorname{Re}\left(e^{-i2\varphi(t)}\left(\sigma^2(t)\right)'\right)=2\bigl((h\cos\varphi-1)(1-t)+t(h^2-h\cos\varphi)\bigr),\\[0.2cm]
\operatorname{Im}\left(e^{-i2\varphi(t)}\left(\sigma^2(t)\right)'\right)=2h\sin\varphi.
\end{array}
\]
On the other hand,
\[
\begin{split}
\bigl(\vert\sigma(t)\vert^2\bigr)'&=(\sigma+\overline{\sigma}-2)(1-t)+t(2\sigma\overline{\sigma}-\sigma-\overline{\sigma})\\
&=2\bigl((h\cos\varphi-1)(1-t)+t(h^2-h\cos\varphi)\bigr).
\end{split}
\]
So the expression of the derivative becomes,
\[
\begin{split}
\kappa'(t)&=\dfrac{2\kappa(0)\Vert\mathbf{v}\Vert^2\Vert \mathbf{w}\Vert ^{3}\vert\mu\vert^2\vert\sigma(t)\vert^{2(2n-4)}}{\Vert T^{n-1}\mathbf{w}
	\Vert^5}\\
&\cdot
\Bigl[(n+1)\bigl(-1+\cos\gamma\operatorname{Im}\left(z(t)\right)\bigr)\bigl((h\cos\varphi-1)(1-t)+t(h^2-h\cos\varphi)\bigr)\Bigr.\\
&\hspace{0.5cm}\Bigl.+3(n-1)h\cos\gamma\sin\varphi \operatorname{Re}\left(z(t)\right)\Bigr].
\end{split}
\] 
where $z(t)=e^{i\theta}e^{i2(n-1)\varphi(t)}$.

\begin{flushright} $\square$ \end{flushright}

We end this section with some examples. The matrix of the next example generates a B\'ezier curve with monotonic curvature for any degree $n\ge 2$.

\begin{figure}[h]
\centering
\begin{subfigure}{0.4\textwidth}
\includegraphics[width=\textwidth]{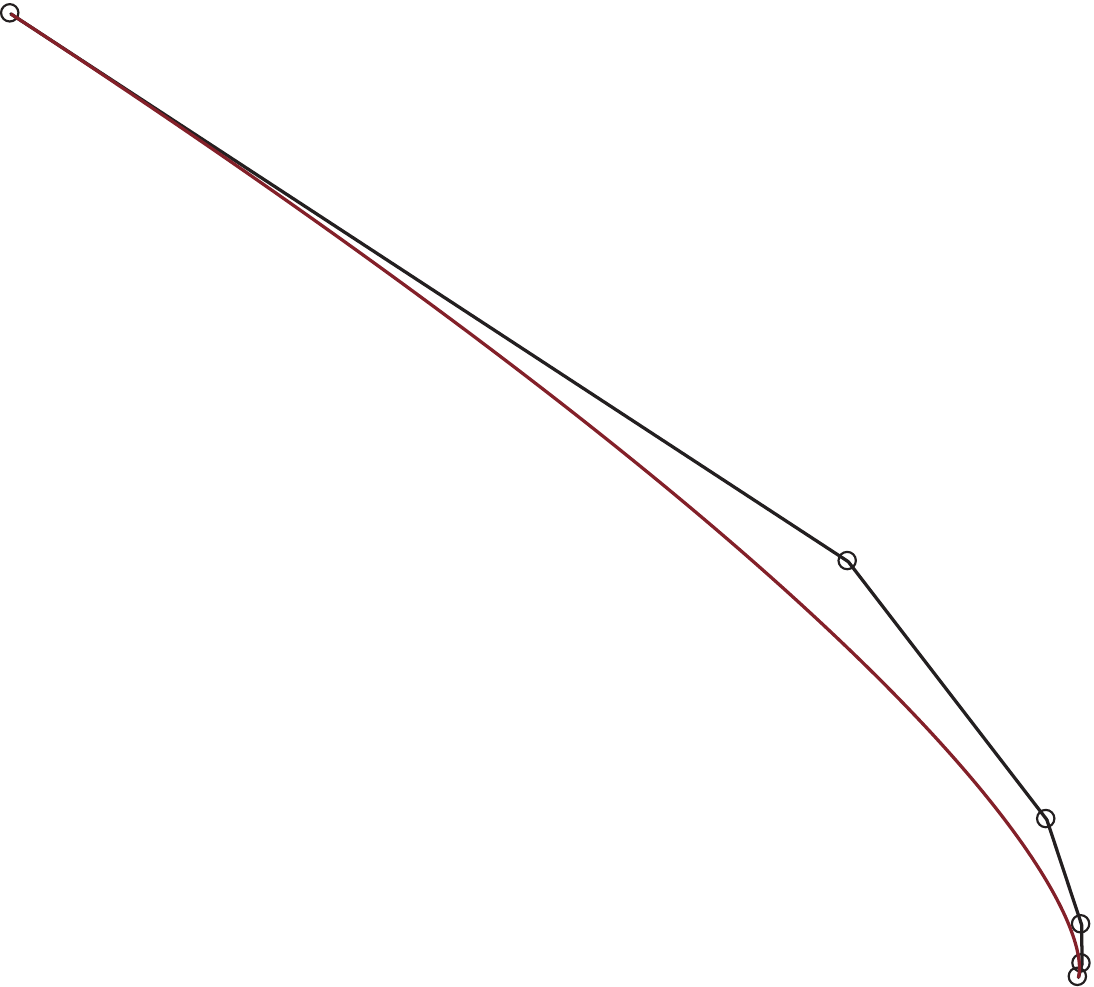}
\end{subfigure}
\hspace{1cm}
\begin{subfigure}{0.3\textwidth}
\includegraphics[width=\textwidth]{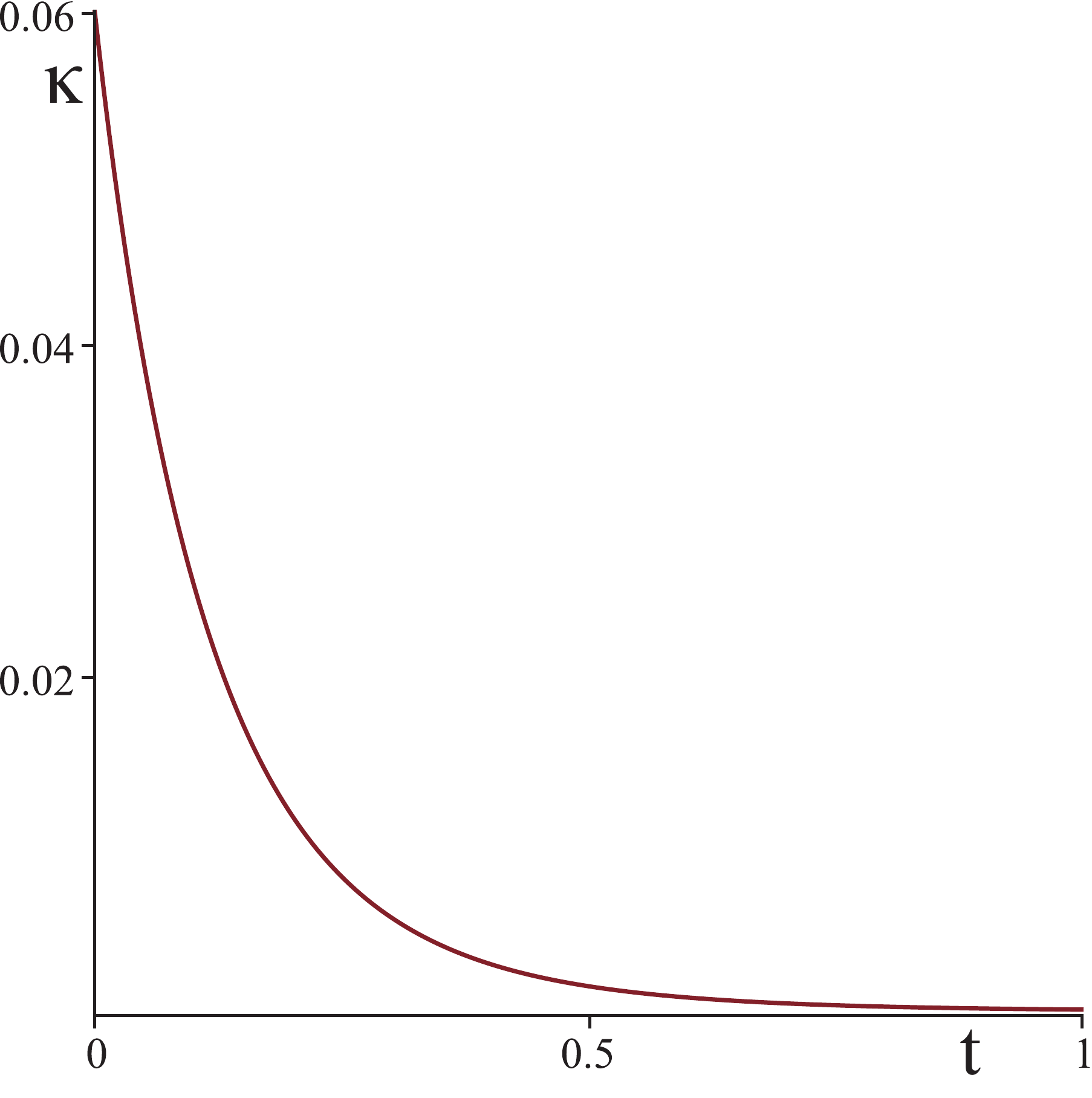}
\end{subfigure}
\caption{Quintic B\'ezier curve generated by $M_{13}$ and 
$\mathbf{w}_{13}$
and the graph of its curvature}\label{eje9}
\end{figure}

\begin{example}
Let 
\[
M_{13}=\frac{3}{4}\begin{pmatrix} 2\sqrt{2+\sqrt{3}}+2-\sqrt{3} &
-\frac{(6-\sqrt{3})\sqrt{2-\sqrt{3}}}{\sqrt{2+\sqrt{3}}} \\[0.3cm] 1 &
2\sqrt{2+\sqrt{3}}-2+\sqrt{3}\end{pmatrix}\approx \begin{pmatrix} 3.1
& -0.86\\ 0.75 & 2.7\end{pmatrix}
\]
with complex eigenvalue $\sigma=3e^{-i\pi/12}$ and $\gamma=5\pi/12$
the angle between the real and the imaginary parts of its eigenvector
$\mathbf{v}= \left(\begin{smallmatrix}1+i\cos\frac{5\pi}{12}\\
i\sin\frac{5\pi}{12}\end{smallmatrix}\right)$.  Thus, condition
(\ref{exMineur}) holds for $M_{13}$.

The B\'ezier curve of degree $5$ generated by $M_{13}$ and
$\mathbf{w}_{13}=10\left(\begin{smallmatrix}\cos(5\pi/12)\\
\sin(5\pi/12)\end{smallmatrix}\right)$, has positive decreasing
curvature.  See Figure \ref{eje9}.

\end{example}

The next example shows a matrix for which the condition in Theorem
\ref{t:complejo} is not valid for any $n\ge 2$.

\begin{figure}[h]
\centering
\begin{subfigure}{0.4\textwidth}
\includegraphics[width=\textwidth]{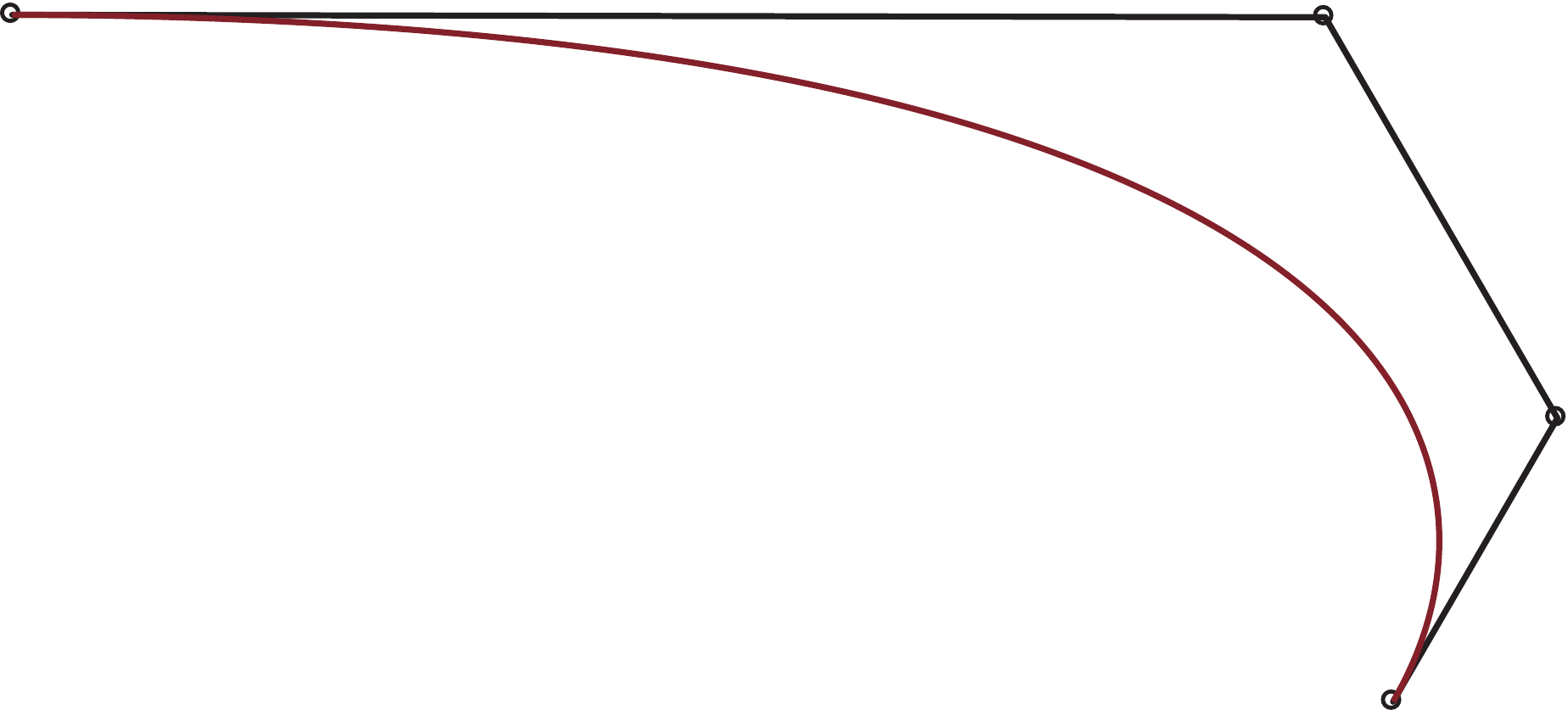}
\end{subfigure}
\hspace{1cm}
\begin{subfigure}{0.3\textwidth}
\includegraphics[width=\textwidth]{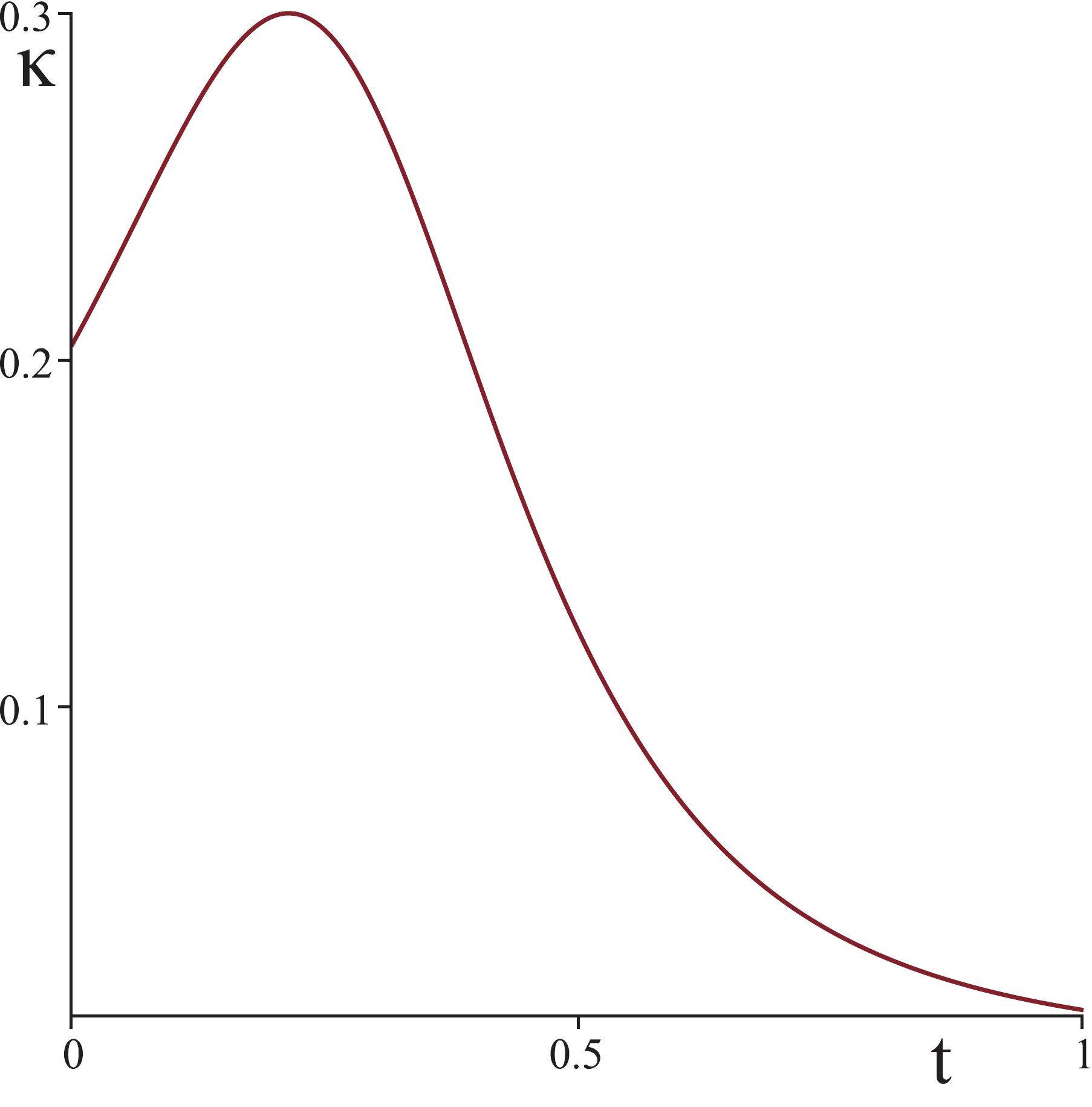}
\end{subfigure}
\caption{Cubic B\'ezier curve generated by $M_{14}$ and 
$\mathbf{w}_{14}$
and the graph of its curvature}\label{eje10}
\end{figure}

\begin{example}
Let
\[
M_{14}=\frac{1}{2}\begin{pmatrix} 3\sqrt{2} &
-5\frac{\sqrt{2}}{\sqrt{3}} \\[0.2cm] \sqrt{6} &
\sqrt{2}\end{pmatrix}\approx \begin{pmatrix} 2.12 & -2.04\\ 1.22 &
0.71\end{pmatrix}
\]
with complex eigenvalue $\sigma=2e^{i\pi/4}$ and $\gamma=\pi/3$ the
angle between the real and the imaginary parts of its eigenvector
$\mathbf{v}= \left(\begin{smallmatrix}1+i\cos\frac{\pi}{3}\\
i\sin\frac{\pi}{3}\end{smallmatrix}\right)$.  Consider
$\mathbf{w}_{14}=\left(\begin{smallmatrix}2\\2\sqrt{3}\end{smallmatrix}\right)$.
The B\'ezier curve of degree $3$ generated by $M_{14}$ and 
$\mathbf{w}_{14}$
does not have monotonic curvature.  See Figure \ref{eje10}.

Notice that condition (\ref{exMineur}) does not hold although
$h\cos\varphi-1=\sqrt{2}-1>0$.

\end{example}

The next example shows a matrix for which the condition
(\ref{exMineur_grado}) holds only for some values of $n$.

\begin{figure}[h]
\centering
\begin{subfigure}{0.15\textwidth}
\includegraphics[width=\textwidth]{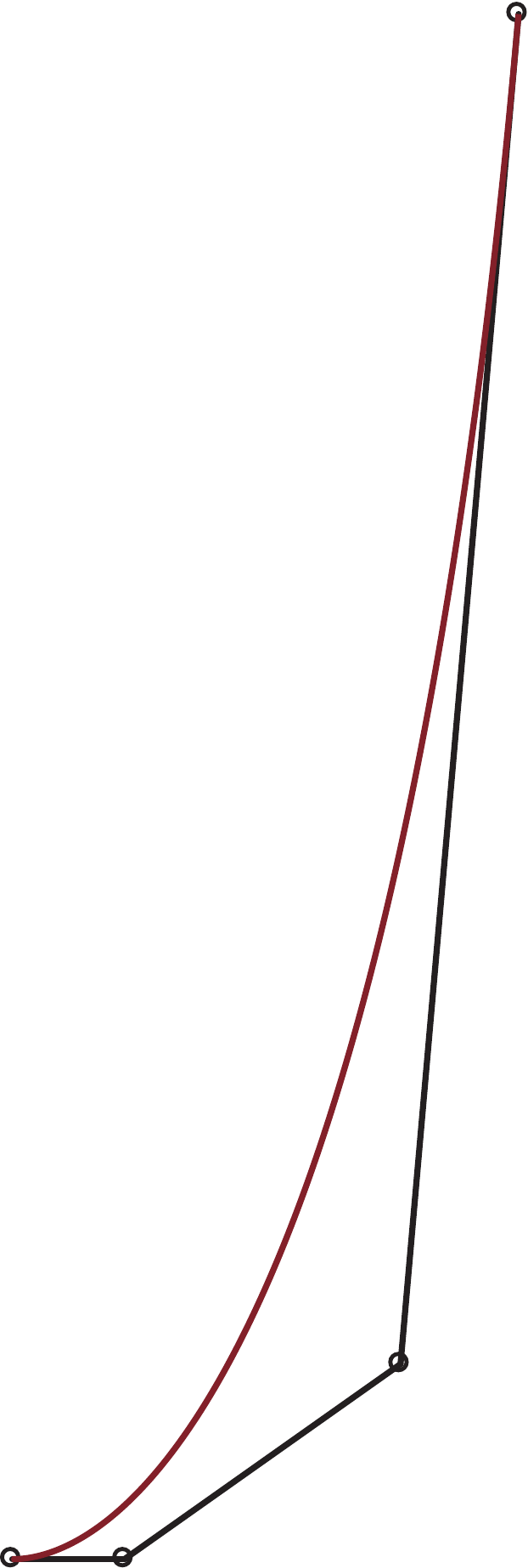}
\end{subfigure}
\hspace{1cm}
\begin{subfigure}{0.3\textwidth}
\includegraphics[width=\textwidth]{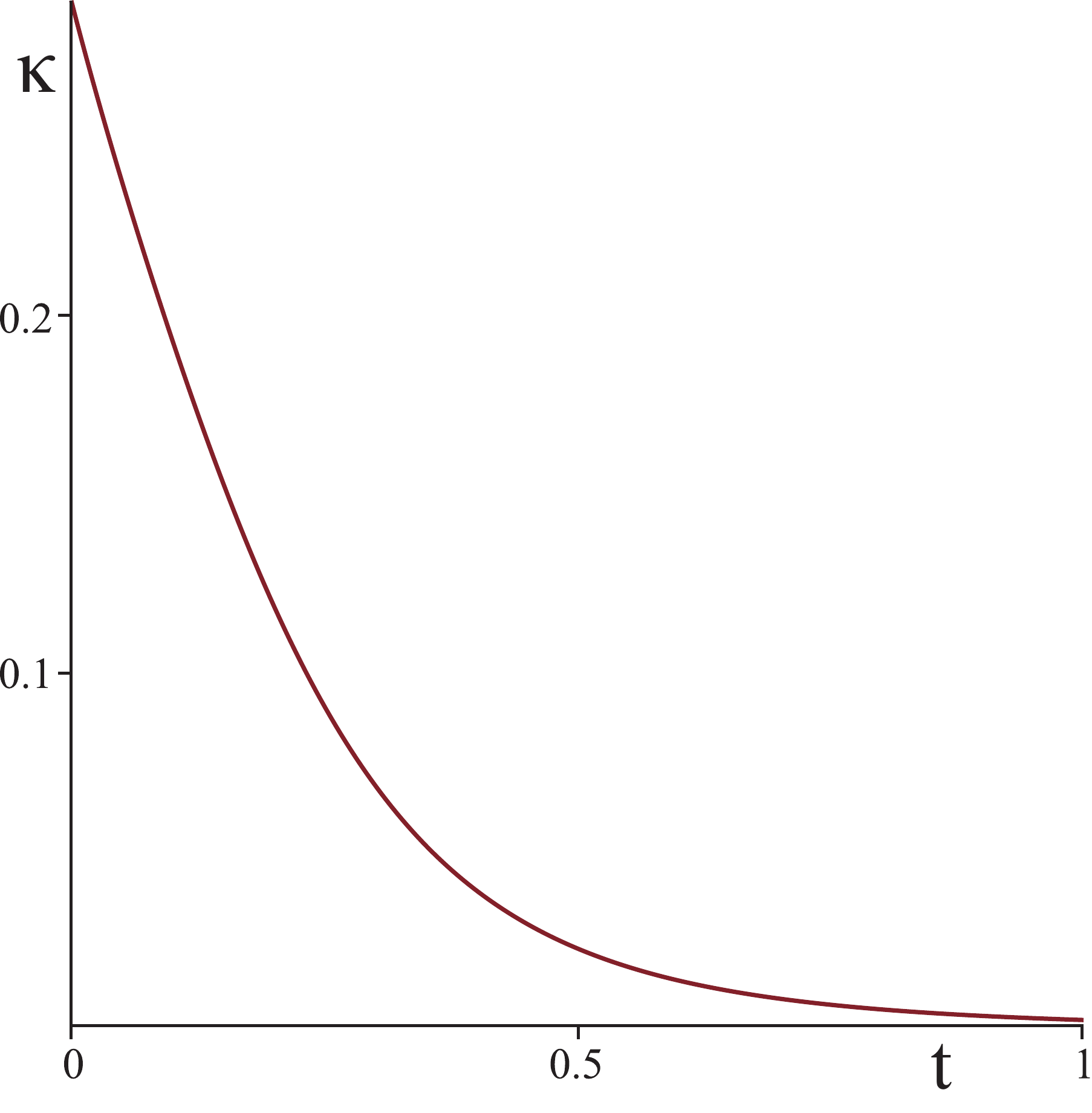}
\end{subfigure}
\caption{Cubic B\'ezier curve generated by $M_{15}$ and 
$\mathbf{w}_{15}$
and the graph of its curvature}\label{eje11}
\end{figure}

\begin{figure}[h]
\centering
\begin{subfigure}{0.5\textwidth}
\includegraphics[width=\textwidth]{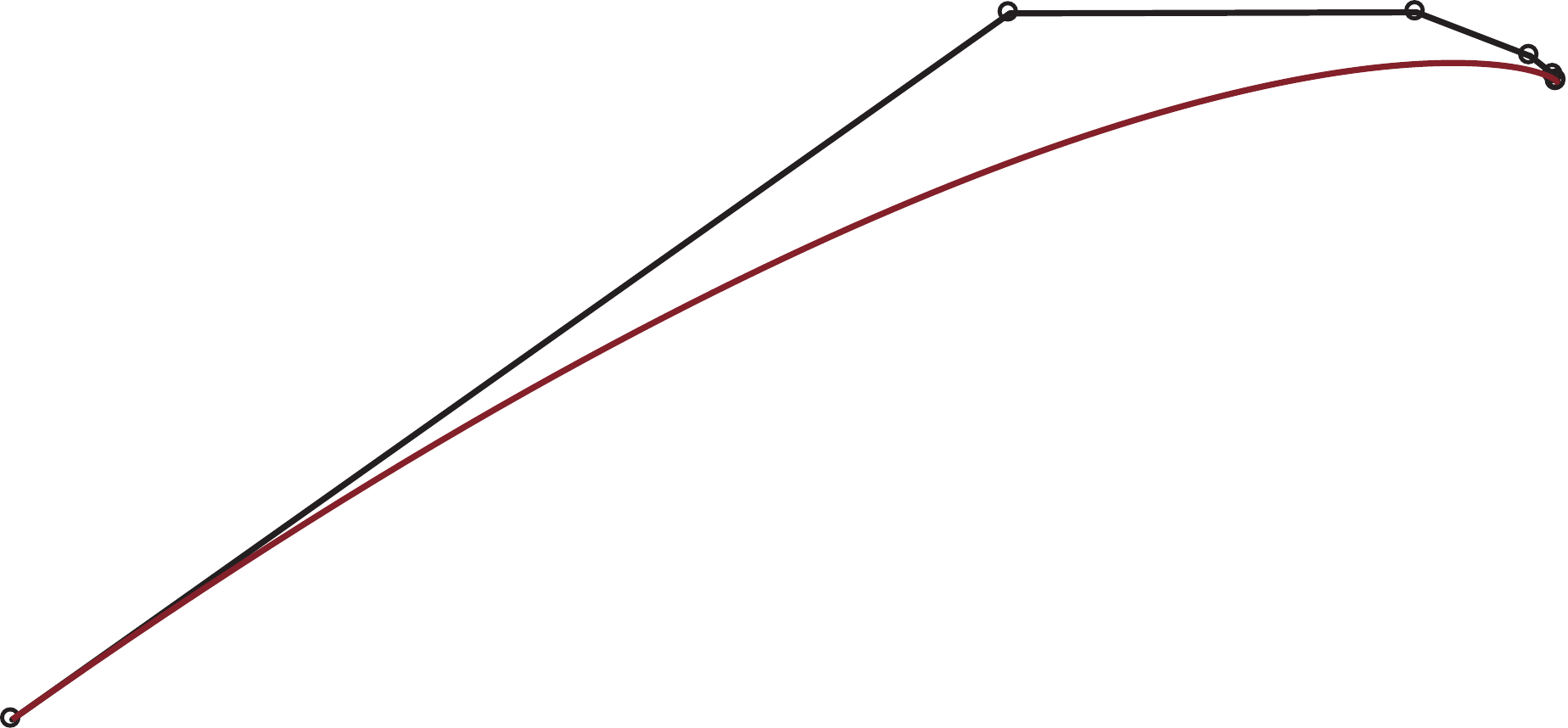}
\end{subfigure}
\hspace{1cm}
\begin{subfigure}{0.3\textwidth}
\includegraphics[width=\textwidth]{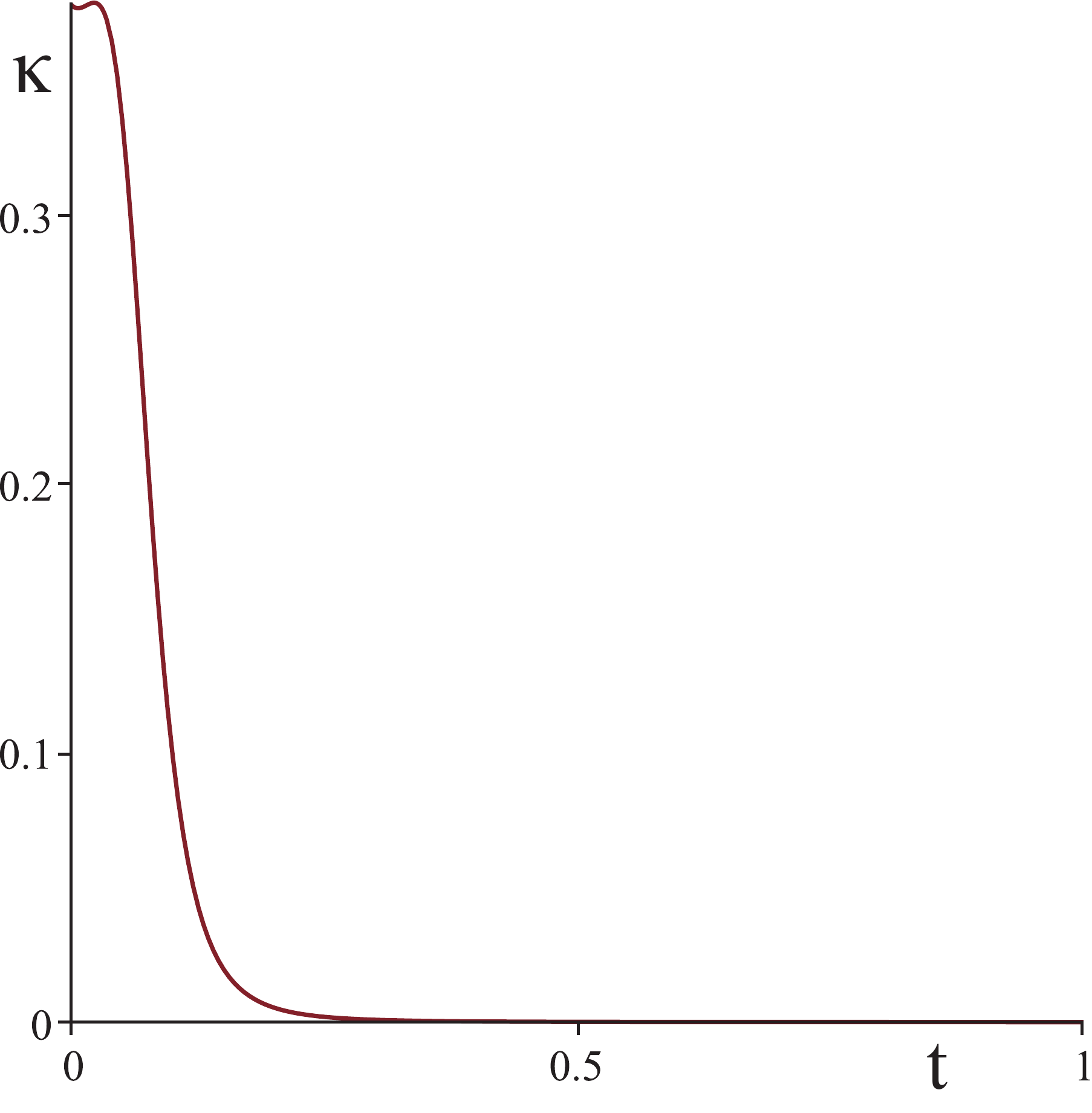}
\end{subfigure}
\caption{B\'ezier curve generated by $M_{15}$ and $\mathbf{w}_{15}$
and the graph of its curvature}\label{eje11a}
\end{figure}

\begin{example}
Let
\[
M_{15}=\begin{pmatrix} 2\sqrt{3}-1 & -\frac{5}{\sqrt{3}} \\[0.2cm]
\sqrt{3} & 2\sqrt{3}+1\end{pmatrix}\approx \begin{pmatrix} 2.46 &
-2.89\\ 1.73 & 4.46\end{pmatrix}
\]
with complex eigenvalue $\sigma=4e^{i\pi/6}$ and $\gamma=2\pi/3$ the
angle between the real and the imaginary parts of its eigenvector
$\mathbf{v}= \left(\begin{smallmatrix}1+i\cos\frac{2\pi}{3}\\
i\sin\frac{2\pi}{3}\end{smallmatrix}\right)$.  Take
$\mathbf{w}_{15}=\left(\begin{smallmatrix}4\\0\end{smallmatrix}\right)$.

The B\'ezier curve of degree $3$ generated by $M_{15}$ and 
$\mathbf{w}_{15}$
has decreasing curvature (see Figure \ref{eje11}) but the B\'ezier
curve of degree $8$ generated by $M_{15}$ and $\mathbf{w}_{15}$ does not
have monotonic curvature (see Figure \ref{eje11a}).  In fact, condition
(\ref{exMineur_grado}) holds for $n=2,\dots,5$ but does not hold for
$n\ge 6$.

\end{example}

\section{Class A B\'ezier curves and related results\label{comparison}}

 In \cite{farin-a} Class A B\'ezier curves are defined as space B\'ezier
 curves for which each edge of the control polygon is obtained from the
 previous one by the action of a Class A matrix, which is a square
 matrix $M$ satisfying certain conditions imposed in order to get curves with monotonic curvature and torsion and, in this way, to generalise the concept of typical curves introduced by Mineur et al. in \cite{mineur}. One of the conditions on $M$ is written in terms of its singular values, that is, the eigenvalues of the symmetric matrix $M^{t}M$. Following Farin´s notation in \cite{farin-a}, in this section $\sigma_1 \geq \sigma_2 \geq \sigma_3$ denote the singular values of the matrix $M$.

In \cite{farin-a} the two conditions that a Class A matrix $M$ is required to satisfy are:

\begin{enumerate}
	
\item For $t \in [0,1]$, and every $\textbf{v}$,
\begin{equation}\label{condition1} \|(1-t)\mathbf{v}+tM\mathbf{v} \|
	\geq \|\mathbf{v} \|,
\end{equation}

(i.e. $M$ is viewed as an ``expansion" matrix).

\item The singular values of the matrix $M$,

\begin{equation}\label{condition2}
\sigma_3^3 \geq \sigma_1,
\end{equation}

(i.e. $M$ is not ``distorting" lengths ``too much").

Note that in \cite[expression (6)]{farin-a} there is a misprint, where there is written $\sigma_3^2 \geq \sigma_1$ instead of $\sigma_3^3 \geq
\sigma_1$, which is what can be deduced from the reasoning in \cite{farin-a}.

\end{enumerate}

We note that condition (\ref{condition1}) is equivalent to $\mathbf{v} \cdot M\mathbf{v }\, \geq \mathbf{v} \cdot \mathbf{v }$. Geometrically, this means that the size of the orthogonal projection of $M\mathbf{v}$ over $\mathbf{v}$ must be at least the size of $\mathbf{v}$, that is:
\begin{equation} \label{farininterpretada}
\| M\mathbf{v} \| \cos \alpha \geq\| \mathbf{v}\|,
\end{equation}
with $\alpha$ being the angle determined by $\mathbf{v}$ and $M\mathbf{v}$.  Written in this way, condition (\ref{condition1}) is a direct generalisation of the condition for decreasing curvature of typical curves of \cite{mineur},  $h \cos \alpha >1$, with the matrix $M$ being the product of a rotation matrix of angle $\alpha$ and a dilatation of factor $h$.

In \cite{farin-a} it is claimed that, with the previous definition, the B\'ezier curves generated by a Class A matrix $M$, that is Class A curves, have two properties that guarantee monotonic curvature and torsion, namely:

\begin{enumerate}
	\item [i)] Class A curves are invariant by subdivision.
	\item [ii)] The curvature at the initial point of a Class A curve is bigger (smaller) than the curvature at the end point. And the torsion behaves in the same way.
\end{enumerate}

However, regarding property i) the inequality (\ref{condition2}),  $\sigma_3^3 \geq \sigma_1$, is not preserved by subdivision of the curve. This is easily seen for a diagonal matrix $M$.  According to
\cite{farin-a}, after subdivision of the interval $[0,1]$ at $t$ the matrix
$M$ is replaced by either $T=(1-t)\mathbb{I}+ Mt$ (to reparameterise the arc of the curve in the interval $[0,t]$) or $MT^{-1}$ (to reparameterise the arc of the curve in the interval $[t,1]$). The
smallest eigenvalue of $T$, $1-t+\sigma_{3}t$, and the largest one,
$1-t+\sigma_{1}t$, are to fulfill (\ref{condition2}). That is,
\[
f(t)=(1-t+\sigma_{3}t)^{3}-(1-t+\sigma_{1}t)
\]
must be non-negative for $t\in[0,1]$ and, since $f(0)=0$, $f'(0)$ must be
non-negative. However,
\[
f'(0)=3\sigma_{3}-\sigma_{1}-2\le 0, \qquad\textrm{ for }
\sigma_{1}\in[3\sigma_{3}-2,\sigma_{3}^{3}],
\]
and hence by choosing
$\sigma_{1}\in[3\sigma_{3}-2,\sigma_{3}^{3}]$, we get counterexamples
of the claim.

\

On the other hand, in \cite{farin-a}, Farin imposes condition (\ref{condition2}) above to show property ii). When considering the curvature, the inequality
$\kappa(0)\ge \kappa(1)$ can be written in terms of the edges of the control polygon and $|M^{j-1}D|$, the area of the triangle formed by $M^{j-1}\mathbf{v}$ and $M^{j}\mathbf{v}$, as
\begin{equation}\label{condition3} 2 \frac{n-1}{n}\frac{| D|}{\| \mathbf{v} \| ^3} \geq 2 \frac{n-1}{n}\frac{|
	M^{n-2} D|}{\| M^{n-1}\mathbf{v} \|^3}.
\end{equation}
%where $\kappa (0)$ and $\kappa (1)$ are expressed in terms of the
%edges of the control polygon %\cite{farin-b}
%and $|M^{j-1}D|$ denotes the area of the triangle formed by $M^{j-1}\mathbf{v}$ and $M^{j}\mathbf{v}$.

The reasoning in \cite{farin-a} begins with the case $n=3$. In this case, proving $\kappa(0)\ge \kappa(1)$ is reduced to show
\[\label{condition4}
\frac{|D|}{\|\mathbf{v} \|^3} \geq \frac{|MD|}{\| M^2\mathbf{v} \| ^3}.
\]

Since $\sigma_1 \geq \sigma_2 \geq \sigma_3$ are singular values of $M$,
\begin{equation}\label{condition5}
\sigma_3 \|\mathbf{v}\| \leq \| M\mathbf{v} \| \leq \sigma_1 \|
\mathbf{v} \|,
\end{equation}
% which is formula $(9)$ in \cite{farin-a}.
and by (\ref{condition1}), $\|M^2\mathbf{v} \| \ge \| M\mathbf{v} \|$, then
\[
\frac{1}{\|\mathbf{v} \|^3} \geq \frac{\sigma_3^3}{\| M^2\mathbf{v} \| ^3}.
\]

The proof in \cite{farin-a} is completed by deducing from (\ref{condition5}) the inequality
\begin{equation}\label{condition10}
\sigma_3 |D|\leq |MD| \leq\sigma_1 |D |,
\end{equation}
which relates the areas $|D|$ and $|MD|$, with the smallest and largest singular values of the matrix $M$ (see \cite[inequality (10)]{farin-a}).
%This claim is supposed to be a consequence of  $(\ref{condition5})$, since $D$ and $MD$ share a common edge.
This would lead to establish condition (\ref{condition2}), $\sigma_{3}^3 \geq \sigma_{1}$, on the singular values of $M$  in order to finish the proof of property ii).

However, it is easy to check that claim (\ref{condition10}) is not valid.  For instance, for typical curves, the matrix $M$ is a rotation with angle $\varphi$ and a scaling with factor $h>1$, hence the singular values are all equal to the factor $h>1$ and the angle between $\mathbf{v}$ and $M\mathbf{v}$ is equal to the angle between $M\mathbf{v}$ and $M^2\mathbf{v}$ and it is $\varphi$.
Hence $|D|=h\|\mathbf{v}\|^{2}\sin\varphi$, $|MD|=h^{3}\|\mathbf{v}\|^{2}\sin\varphi$,
and (\ref{condition10}) would imply
% \begin{equation}
% \sigma_3 h \, \leq \, h^3 \, \leq \,  \sigma_1 h
% \end{equation}
\[
h^2  \leq h^3 \leq h^2 ,
\]which is impossible  for $h>1$.

In fact, using the bounds for $\|M\mathbf{v} \|$ in terms of the singular values of $M$ given by
(\ref{condition5}),
\begin{equation}\label{acotacionM2v}
\sigma_3^m \|\mathbf{v}\|\leq \|M^m\mathbf{v} \| \leq \sigma_1^m \| \mathbf{v} \|,
\end{equation}
and therefore a valid expression relating the areas $|D|$ and $|MD|$, (instead of \cite[inequality (10)]{farin-a}), is
\begin{equation}\label{condition10correcta} \sigma_3^2 \,| D |\,
\frac{\sin \alpha_2}{\sin \alpha_1}  \leq | MD | \leq \sigma_1^2 \,
|D|\, \frac{\sin \alpha_2}{\sin \alpha_1},
\end{equation}
where $\alpha_j$ is the angle determined by the  vectors $M^{j-1}\mathbf{v}$ and
$M^{j}\mathbf{v}$.

Now, following the proof of \cite[inequality (8)]{farin-a} but using (\ref{acotacionM2v}) with $m=2$ and (\ref{condition10correcta}), and extending to arbitrary degrees by iteration of the previous bounds, we get:

\begin{proposition}\label{condicion2correcta} If $M$ satisfies condition $(\ref{condition1})$ and its singular values satisfy
	\[
\sigma_3^{3(n-1)} \geq \sigma_1^{2(n-2)} \, \frac{\sin \alpha_{n-1}}{\sin \alpha_1}
	\]
	with $\alpha_j$ the angle between $M^{j-1}\mathbf{v}$ and $M^{j}\mathbf{v}$ then,
	\[
	\frac{|D|}{\| \mathbf{v} \|^3} \geq \frac{| M^{n-2}D |}{\| M^{n-1}\mathbf{v}
		\|^3}\]
	where $| M^{j-1}D |$ is the area of the triangle determined by $M^{j-1}\mathbf{v}$ and $M^{j}\mathbf{v}$.
\end{proposition}

In particular, for $n=3$, if the singular values of the matrix $M$ satisfy
\[
\sigma_3^6 \geq \sigma_1^2 \, \frac{\sin \alpha_2}{\sin \alpha_1}
\]
(instead of $\sigma_3^3 \geq \sigma_1$) then
\[
\frac{|D|}{\| \mathbf{v} \|^3} \geq \frac{| MD |}{\| M^2\mathbf{v}\|^3},
\]
as it was needed in \cite{farin-a} to show that $\kappa (0) \geq \kappa (1)$, where $\kappa(t)$ is the curvature of a cubic B\'ezier curve whose edges of the control polygon are given by $M^{j}\mathbf{v}$ for $j=0,1,2$.

\

We finish this section with some comments about some results in the literature regarding Farin's Class A curves.

\

In \cite{cao} an example of a Class A matrix is provided for which the condition involving the singular values of the
matrix, $\sigma_{1}=1.102$ and $\sigma_{3}=1.05$, is not supposed to be invariant under subdivision at $t=0.5$. However this is a counterexample for the misprinted condition $\sigma_{3}^{2}\ge \sigma_{1}$, but not for $\sigma_{3}^{3}\ge \sigma_{1}$.

\

% Farin2006, used by Farin to prove that for a Class A curve, the
% curvature at initial point is at least the curvature at the final
% point.

% Hence it is not explicitly shown that the constructed Class A
% B\'ezier curves with them have non-monotonic curvature plot.

In \cite{zhao} two counterexamples are found for conditions
(\ref{condition1}) and (\ref{condition2}) on the Class A matrix that
generate a Class A B\'ezier curve.  Note that the matrix $M$ is 
$M^{t}$ in our notation.  The first example of a cubic B\'ezier curve is
not a counterexample for Farin's method because it does not satisfy
condition (\ref{condition1}) since the given
$\mathbf{v}=\begin{pmatrix}0.9724\\ 0.2333\end{pmatrix}$ and $M =
\begin{pmatrix}
1.2545 & -2.9594 \\ 
1.5576 & 2.3836
\end{pmatrix}$
produce
\[
\frac{\mathbf{v} \cdot M\mathbf{v}}{\mathbf{v} \cdot \mathbf{v}}
=0.9979<1.
\]

\

In the second example of \cite{zhao} the curve has non-monotonic
curvature although it satisfies the misprinted condition in
\cite{farin-a} about singular values of $M$ and it also satisfies
Proposition \ref{condicion2correcta} here, showing that Farin's method
does not work.  This points out that the condition about singular
values involved in Proposition \ref{condicion2correcta} is not
preserved under subdivision of the curve.

\section*{Acknowledgments}

This work is partially supported by the Spanish Ministerio de
Econom\'\i a y Competitividad through research grant TRA2015-67788-P.

\bibliographystyle{elsarticle-num}
\bibliography{cagd}

\end{document}